\documentclass[12pt,twoside]{article}
\usepackage{amsmath,amssymb,amsthm}
\usepackage{graphicx}
\usepackage{mathrsfs}
\usepackage{bm}
\usepackage{color}
\usepackage{enumerate}
\usepackage{enumitem}
\usepackage{upgreek}
\usepackage{verbatim}
\usepackage[sort&compress,numbers]{natbib}
\makeatletter
\renewcommand{\@makefnmark}{\hbox{}}
\renewcommand{\@makefntext}[1]{\noindent\reset@font\footnotesize#1}
\makeatother

\setlist[description,1]{font=\textup,
	leftmargin=7mm,labelsep=1.5mm,itemsep=0.2mm,parsep=0.2mm}
\setlist[description,2]{font=\textup,
	leftmargin=7mm,labelsep=1.5mm,topsep=0.2mm,itemsep=0.2mm,parsep=0.2mm}

\newtheorem{definition}{\noindent\textbf{Definition}}[section]
\newtheorem{lemma}{\noindent\textbf{Lemma}}[section]
\newtheorem{theorem}{\noindent\textbf{Theorem}}[section]
\newtheorem{remark}{\noindent\textbf{Remark}}[section]

\numberwithin{equation}{section}

\allowdisplaybreaks
\makeatletter
\evensidemargin
\oddsidemargin
\makeatother

\title{\Large\bf  A new partial differential nonlinear system containing quasivariational and parabolic variational inequalities and its application}

\author{Wei Li$^{a,b}$, \ Zhenghui Tang$^{a,c}$, \ Zengbao Wu$^{d}$\footnote{Corresponding author. E-mail address: zengbaowu@hotmail.com,  wuzengbao@lynu.edu.cn}, \ Chunyan Yang$^{e}$  \\
	$^{a}${\small \textit{School of Mathematical Sciences, Chengdu University of Technology,}}\\
	{\small\textit{Chengdu, Sichuan 610059, P.R. China}}\\
	$^{b}${\small \textit{Mathematics Research Center, Chengdu University of Technology, }}\\
	{\small\textit{Chengdu, Sichuan 610059, P.R. China}}\\
	$^{c}${\small \textit{Geomathematics Key Laboratory of Sichuan Province, Chengdu University of Technology, }}\\
	{\small\textit{Chengdu, Sichuan 610059, P.R. China}}\\
	$^{d}${\small \textit{Department of Mathematics, Luoyang Normal University, }}\\
	{\small\textit{Luoyang, Henan 471934, P.R. China}}\\
$^{e}${\small \textit{Department of Mathematics, Sichuan University, Chengdu, Sichuan 610064, P.R. China}}
}
\date{ }

\begin{document}
\maketitle

\begin{abstract}
\noindent 
We study a new nonlinear system which contains a partial differential equation, 
a quasivariational inequality and a parabolic variational inequality in Banach spaces. 
We obtain the unique solvability of the coupled system under moderate conditions by using the Banach's fixed point theorem. We employ the main results to investigate a viscoelastic frictional contact problem with long-memory effects, wear processes, and damage phenomenon.
\end{abstract}

\footnote{\,This work was supported by National Natural Science Foundation of China (12301395, 11901273), the Program for Science and Technology Innovation Talents in Universities of Henan Province (23HASTIT031), and the Natural Science Foundation of Henan (252300421997).}
\vskip 0.2cm

\section{Introduction}
	Variational inequalities (VIs) have been widely used in engineering, optimization, economics, mechanics and other scientific fields. When ordinary differential equations (ODEs) are introduced, the system consisting of VIs and ODEs is known as differential variational inequalities (DVIs) and this concept was introduced and studied in depth by Pang et al. in~\cite{pang2008differential}. Then, Chen in~\cite{chen2013convergence} used a numerical method to solve a class of DVIs and proved the convergence of this method. Furthermore, Liu et al. in~\cite{liu2017partial} focused on partial differential variational inequalities (PDVIs) and proved the solvability to this class of problems utilizing the semigroup theory and fixed point theorems. Zhang et al. in~\cite{zhang2023penalty} treated the solvability of stochastic DVIs. Hao et al. in~\cite{hao2025existence} analyzed the unique solvability and convergence of the solution to a class of generalized second-order delay differential variational-hemivariational inequalities, with relevant applications. Mig\'{o}rski in~\cite{migorski2021} studied an abstract system of two evolution inclusions which could be transformed into two classes of variational inequalities and provided important support for scholars researching on multiple quasistatic or evolutionary inclusions. Currently, DVIs and PDVIs have aroused widespread interest in relevant academia owing to their crucial applications in frictional contact models,  dynamic transportation, network optimization and Nash equilibrium, and significant results have been made at both the theoretical and applied levels in \cite{liu2016evolutionary,pang2009solution,liu2019differential,Liu2018,migorski2018hemi,2025CT,Li2020,Li2015,ceng2025fractional,ceng2025evolution}.
	
	On the other hand, wear phenomena in frictional contact problems have become a topic of great interest in the industrial field. Andrews et al. in~\cite{andrews1997dynamic} studied the frictional contact problem involving wear. Furthermore, they modeled the dynamic thermo-viscoelastic contact problem by Archard's wear law and proved the unique solvability of the problem. In recent years, Sofonea et al. in~\cite{sofonea2016analysis} developed a model for the quasi-static elastic contact problem based on the generalized form being used in \cite{andrews1997dynamic}, and further proved the unique solvability of the problem. Chen et al. in~\cite{chen2020variational} modeled the contact problem by Archard's wear law, and then obtained the unique solvability of the problem. However, the generalized differential Archard condition allows the wear diffusion, which was studied in planar contact surface. So, in a frictional contact problem, the total wear can be better described by partial differential equations, see \cite{kalita2019frictional}.
	As the load-bearing capacity of the materials, such as cement and rubber, decreases over time, the phenomenon is known as damage of materials, see \cite{fremond1995damage}. The damage function was used for the first time in \cite{fraemond1996damage} to quantify the damage of the materials. Many interesting models related to the damaged contact problems can be found in \cite{han2020numerical,li2010quasistatic,gasinski2015variational,szafraniec2016dynamic,sofonea2005analysis,xuan2021numerical}. In recent years, Gasi\'{n}ski et al. in~\cite{gasinski2015variational} studied a contact problem involving damage using the method of variational-hemivariational inequalities.
	
	In frictional contact problems, materials with long memory effects, including rubbers and pastes, are transformed to problems with the history-dependent operators. Some examples and explanations can be found in \cite{drozdov1996finite,banks2011brief}. A large number of frictional or frictionless models for quasi-static contact have given rise to many history-dependent VIs and evolutionary VIs, and the models have already been investigated in some abstract frameworks. Mig\'{o}rski in~\cite{migorski2022} studied a class of variational-hemivariational inequality coupled with an ODE, and applied it to the dynamic unilateral viscoplastic frictionless contact problem and the viscoelastic contact problem with friction and adhesion. Sofonea et al. in~\cite{sofonea2011history} studied a general class of quasivariational inequalities emerging from numerous mathematical models describing quasistatic contact processes. The results on  quasivariational inequalities from \cite{sofonea2011history} have been applied in the research of a large number of quasi-static contact models. Wang et al. in~\cite{wang2025class} studied the solvability of a type of time-dependent mixed quasivariational-hemivariational inequality and applied the results to a contact model for elastic materials and the Oseen model of generalized incompressible Newtonian fluids. Liang et al. in~\cite{liang2024second} studied a second-order differential inclusion driven by a quasivariational-hemivariational inequality with a perturbation operator in Banach spaces.
	
	The works mentioned above primarily focus  on double-coupled systems that incorporate a VI. This VI couples a differential equation (such as an ordinary differential equation, partial differential equation, or fractional differential equation), where the VI could be quasivariational inequality or hemivariational inequality, etc. Furthermore, Chen et al. in~\cite{CT2021} studied a class of three-coupled ordinary differential systems involving a ordinary differential equation, a variational-hemivariational inequalities and a parabolic inequality, and applied it to the plane contact problem. Inspired by these papers, we investigate a class of three-coupled partial differential systems containing a partial differential equation, a quasi-variational inequality and a parabolic inequality, and apply it to the curve contact problem. 
	
	Let $\varphi\colon V\times Y\to\mathbb R$, $A\colon X\times K_V\to V$, $\mathcal R\colon C(I;K_V)\times C(I;K_Y)\to C(I;X)$, $\mathcal S\colon C(I;K_V)\to C(I;Z)$, $j\colon Z\times V\times\mathbb R\times K_V\to \mathbb R$, $f\colon I\to V$, $g\colon Y\times Y\to\mathbb R$, $\phi\colon I\times V\to Y_1$. 
	The problem under investigation reads as folows: 
	Find $\eta\colon I\to K_V$, $\xi \colon I\to K_Y$ and $w\colon I\times\Omega\to\mathbb R$ satisfying the following inequalities and equations
	\begin{align}\label{0}
		\frac{\partial w(t,x)}{\partial t}-\Delta 	w(t,x)=\varphi(\eta(t),\xi(t)),\\
		\begin{aligned}\label{1}
			\langle A(\mathcal R(&\eta(t),\xi(t)),\eta(t)),v-\eta(t)\rangle_V+j(\mathcal S\eta(t),\eta(t),w(t,x),v)\\[1mm]
			&-j(\mathcal S\eta(t),\eta(t),w(t,x),\eta(t))\geq\langle f(t),v-\eta(t)\rangle_V \text{ for all $v\in K_V$},
		\end{aligned}\\		\label{2}		\langle\dot{\xi}(t),\delta-\xi(t)\rangle_{Y_1}+g(\xi(t),\delta-\xi(t))\geq\langle \phi(t,\eta(t)),\delta-\xi(t)\rangle_{Y_1}\text{ for all $\delta\in K_Y$},\\
		\frac{\partial w(t,x)}{\partial\nu}+bw(t,x)=0 \text{ for all $x\in\partial\Omega$},\\
		w(0,x)=w_0(x), \ \ \xi(0)=\xi_0,
	\end{align}
	for each $t\in I$ and $x\in\Omega\subset\mathbb R^n$, $n=2$, $3$, where $b\in\mathbb{R}$, $X,\ Y,\ Z$ are reflexive, separable Banach spaces, and $X^*,\ Y^*,\ Z^*$ are the dual spaces of $X,\ Y,\ Z$, respectively. Moreover, we  suppose that $Y_1$ and $V$ are separable Hilbert spaces satisfying the embedding $Y\subset Y_1\subset Y^*$, the sets $K_V\subset V$ and $K_Y\subset Y$ are convex, and $I=[0,T]$ with $T\in\mathbb{R^+}$.

In the paper, our main novel contributions are as follows.
(i) We investigate a new class of nonlinear systems which contains a partial differential equation, a quasivariational inequality and a parabolic variational inequality in Banach spaces, which is distinct from, and extends  previous double-coupled systems, see, e.g., \cite{migorski2021,migorski2022}. 
(ii) Most researchers use the operator semigroup theory to study partial differential equations and use the fixed point theorem to solve the class of system, see, e.g., \cite{liu2017partial,Liu2018,Li2020}, 
while we study the partial differential equation in this system by employing the estimates on solutions and use the Banach's  fixed point theorem to prove unique solvability of the system. 
(iii) We apply the abstract results for the system to solve a frictional contact problem and establish its unique solvability, where the contact surface is curved and distinct from planar surface of the most studies mentioned above, see, e.g., \cite{liu2019differential,sofonea2016analysis,CT2021}.
	
The structure of this paper is as follows. In Section~\ref{S2}, we recall key definitions and useful results. In Section~\ref{S3}, we obtain the unique solvability of the problem (1.1)-(1.5). In Section~\ref{S4}, we study a specific frictional contact problem, and this problem is formulated as the system (1.1)-(1.5). Furthermore, by applying the theoretical findings derived in  Section~\ref{S3}, we prove the unique solvability to the viscoelastic frictional contact problem with long-memory effects, wear processes, and damage.

\section{Preliminaries}\label{S2}

\begin{definition}
We assume that $\mathcal X$ and $\mathcal Y$ are nonempty sets. The graph, $Gr(S)$, of a single valued map $S\colon\mathcal X\to\mathcal Y$ is defined as follows
$$
Gr(S)=\{(a,b)\in\mathcal X\times\mathcal Y \mid b=Sa\}.
$$
\end{definition}
		
\begin{lemma}\cite[p.9]{sofonea2017variational}
Let $\mathcal X$ and $\mathcal Y$ 
be nonempty sets, 
$\mathcal P\subset\mathcal X\times\mathcal Y$, and 
$Ab=\{a\in\mathcal X \mid (a,b)\in\mathcal P\}$ for all $b\in\mathcal Y$. 
We assume $S\colon \mathcal X\to\mathcal Y$ is a single valued map, and for each $b \in \mathcal{Y}$, there is a unique $a \in \mathcal{X}$ such that the pair $(a, b)\in\mathcal{P}$. 
Then $P\cap Gr(S)\neq\emptyset$ is a singleton if and only if $\Lambda := S A$ 
has a unique fixed point.
\end{lemma}
	
\begin{lemma}\cite[p.21]{sofonea2017variational}
We suppose $X$ is a Hilbert space, $y \in X$, the nonempty closed convex set $K\subset X$, the operator $A \colon K\to X$ is strongly monotone and Lipschitz continuous, and the functional $j\colon K \to \mathbb{R}$ is lower semicontinuous. 
Then there is a unique $b$ satisfying
\begin{equation*}
b\in K,\ \langle Ab,a-b\rangle_X+j(a)-j(b)\geq\langle 	y,a-b\rangle_X \ \ \text{for all} \ \ a\in K.
\end{equation*}
\end{lemma}

\begin{definition}
We suppose $I$ is a time interval and the almost history-dependent operator $S\colon C(I;X)\to C(I;Y)$ is defined as follows: 
for each compact $K\subset I$, there are $l_K\in[0,1)$ and $L_K>0$ such that 
for each $t\in K$,
\begin{equation*}
\|Sa(t)-Sb(t)\|_Y\leq l_K\|a(t)-b(t)\|_X
+ L_K \int_0^t\|a(s)-b(s)\|_X \,ds
\end{equation*}
for all $a$, $b\in C(I;X)$.
\end{definition}
	
\begin{lemma}\cite[p.41]{sofonea2017variational}
We suppose $X$ is a Banach space, $I$ is a time interval and the operator $A\colon C(I;X)\to C(I;X)$ is almost history-dependent. Then $A$ has a unique fixed point.
\end{lemma}

\begin{lemma}\cite[p.8]{nittka2014inhomogeneous}
We suppose $\Omega\subset\mathbb R^N$ is a bounded connected domain with the boundary $\partial\Omega$, $N\geq 2$. 
Let $a_{ij}\in L^\infty(\Omega)$, 
$b_j$, $c_i\in L^q(\Omega)$, $d\in L^{\frac{q}{2}}(\Omega)$ and $\beta\in L^{q-1}(\partial\Omega)$ be given, where $q>N$ is arbitrary. 
We assume there exists $m>0$ such that $\sum_{i,j=1}^{N}a_{ij}\zeta_i\zeta_j\geq m|\zeta|^2$ for all $\zeta\in\mathbb R^N$, 
and 
$T>0$, $w_0\in L^2(\Omega)$, $h\in L^2([0,T];L^2(\Omega))$ and $p\in L^2([0,T];L^2(\partial\Omega))$. 
Then there is a unique weak solution $w\in C([0,T];L^2(\Omega))\cap L^2([0,T];H^1(\Omega))$ of the following Robin problem.
\begin{equation}\label{2.1}
\left\{
\begin{aligned}
w_t(t,x)-Lw(t,x)&=h(t,x),\ t>0,x\in\Omega,\\
\frac{\partial w(t,x)}{\partial\nu_L}+\beta w(t,x)&=p(t,x),\ 	t>0,x\in\partial\Omega,\\
w(0,x)&=w_0(x),\ x\in\Omega,
\end{aligned}
\right.
\end{equation}
where $\nu=(\nu_j)_{j=1}^N$ denotes the outer unit normal on the boundary $\partial\Omega$ 
and
\begin{gather*}						Lw:=\sum_{j=1}^{N}D_j\left(\sum_{i=1}^{N}a_{ij}D_iw+b_jw\right) -\left(\sum_{i=1}^{N}c_iD_iw+dw\right),\\
\frac{\partial w}{\partial\nu_L}:= \sum_{j=1}^{N}\left(\sum_{i=1}^{N}a_{ij}D_iw+b_jw\right)\nu_j, 
\end{gather*}	
where $D_iw=\frac{\partial w}{\partial x_i}$.
\end{lemma}

\begin{lemma}\cite[p.9]{nittka2014inhomogeneous}
We suppose that $T>0$, $a_1$, $a_2$, $b_1$, $b_2\in[2,+\infty)$ 
satisfy $\frac{1}{a_1}+\frac{N}{2b_1}<1$ and $\frac{1}{a_2}+\frac{N-1}{2b_2}<\frac{1}{2}$,   $h\in L^{a_1}([0,T];L^{b_1}(\Omega))$, 
$p\in L^{a_2}([0,T];L^{b_2}(\partial\Omega))$ and $w_0\in L^\infty(\Omega)$. Then the weak solution $w$ of (\ref{2.1}) satisfies
\begin{equation*}
\|w\|_{L^\infty([0,T];L^\infty(\Omega))}\leq c(\|w_0\|_{L^\infty(\Omega)} +\|h\|_{L^{a_1}([0,T];L^{b_1}(\Omega))}+\|p\|_{L^{a_2}([0,T];L^{b_2}(\Omega))}),
\end{equation*}
where $c$ depends only on $T$, $N$, $\Omega$, $a_1$, $b_1$, $a_2$, $b_2$ and the coefficients of the equation.
\end{lemma}

\begin{lemma}\cite[p.10]{nittka2014inhomogeneous}
We assume that $T>0$, $h$ and $p$ satisfy the conditions of Lemma 2.5 and $w_0\in C(\overline{\Omega})$. 
Then the weak solution $w$ of (\ref{2.1}) is in $C([0,T];C(\overline{\Omega}))$. Moreover,  $w(t)\to w_0$ uniformly on $\overline{\Omega}$ as $t\to 0$.
\end{lemma}

\section{Unique solvability for the new nonlinear system}\label{S3}

In this section we prove the unique solvability of the problem (1.1)-(1.5). We need the following assumptions.
	\begin{description}
		\item[\normalfont H($A$)]: $A\colon X\times K_V\to V$ satisfies
		\begin{description}
			\item[\normalfont (a)] $A(\cdot,\cdot)$ is Lipschitz continuous with $L_1$ and $L_2$, i.e., there exist $L_1>0,\ L_2>0$ satisfying
				$$\|A(a_1,b)-A(a_2,b)\|_V\leq L_1\|a_1-a_2\|_X$$
			and
				$$\|A(a,b_1)-A(a,b_2)\|_V\leq L_2\|b_1-b_2\|_V$$
			for all $a_1,\ a_2,\ a\in X$ and $b_1,\ b_2,\ b\in K_V$;
		\end{description}
		\begin{description}
			\item[\normalfont (b)] there exists $m>0$ satisfying
				$$\langle A(a,b_1)-A(a,b_2),b_1-b_2\rangle_V\geq m\|b_1-b_2\|_V^2$$
			for all $a\in X$ and $b_1,\ b_2\in K_V$.
		\end{description}
	\end{description}

	\begin{description}
		\item[\normalfont H($j$)]: $j\colon Z\times V\times\mathbb R\times K_V\to\mathbb R$ satisfies
		\begin{description}
			\item[\normalfont (a)] for all $a\in Z$, $b\in V$ and $c\in\mathbb R$, $j(a,b,c,\cdot)$ is convex and lower semicontinuous on $K_V$;
			\item[\normalfont (b)] there exist $\alpha>0,\ \beta>0,\ \gamma>0$ satisfying
		\begin{eqnarray*}
			&& j(a_1,b_1,c_1,d_2)-j(a_1,b_1,c_1,d_1)+j(a_2,b_2,c_2,d_1)-j(a_2,b_2,c_2,d_2) \nonumber\\
			&\leq&\alpha\|a_1-a_2\|_Z\|d_1-d_2\|_V+ \beta\|b_1-b_2\|_V\|d_1-d_2\|_V+\gamma|c_1-c_2|\|d_1-d_2\|_V
		\end{eqnarray*}
		for all $a_1,\ a_2\in Z,\ b_1,\ b_2\in V,\ c_1,\ c_2\in\mathbb R$ and $d_1,\ d_2\in K_V$.
		\end{description}
	\end{description}

	\begin{description}
		\item[\normalfont H($\mathcal R$)]: $\mathcal R\colon C(I;K_V)\times C(I;K_Y)\to C(I;X)$ satisfies that for any compact set $J\subset I$, there exist $r_{1J}>0,r_{2J}>0$ satisfying
			$$\|\mathcal R(a_1(t),b(t))-\mathcal R(a_2(t),b(t))\|_X\leq r_{1J}\int_0^t\|a_1(s)-a_2(s)\|_Vds$$
		and
			$$\|\mathcal R(a(t),b_1(t))-\mathcal R(a(t),b_2(t))\|_X\leq r_{2J}\int_0^t\|b_1(s)-b_2(s)\|_{Y_1}ds$$
		for all $a_1,\ a_2,\ a\in C(I;K_V),\ b_1,\ b_2,\ b\in C(I;K_Y)$ and $t\in J$.
	\end{description}

	\begin{description}
		\item[\normalfont H($\mathcal S$)]: $S\colon C(I;K_V)\to C(I;Z)$ satisfies that for any compact set $J\subset I$, there exists $s_J>0$ satisfying
			$$\|\mathcal Sa_1(t)-\mathcal Sa_2(t)\|_Z\leq s_J\int_0^t\|a_1(s)-a_2(s)\|_Vds$$
		for all $a_1,\ a_2\in C(I;K_V)$ and $t\in J$.
	\end{description}

	\begin{description}
		\item[\normalfont H($f$)]: $f\colon I\to V$ is a continuous function.
	\end{description}

	\begin{description}
		\item[\normalfont H($\varphi$)]: $\varphi\colon V\times Y\to\mathbb R$ satisfies
		\begin{description}
			\item[\normalfont (a)] $\varphi(\cdot,\cdot)$ is Lipschitz continuous with $L_\varphi$, i.e., there exists $L_\varphi>0$ satisfying
				$$|\varphi(a_1,b_1)-\varphi(a_2,b_2)|\leq L_\varphi(\|a_1-a_2\|_V+\|b_1-b_2\|_{Y_1})$$
			for all $a_1,\ a_2\in V$ and $b_1,\ b_2\in Y\subset Y_1$;
		\end{description}
		\begin{description}
  			 \item[\normalfont (b)] $\varphi(0_V,0_Y)<+\infty$.
		\end{description}
	\end{description}

	\begin{description}
		\item[\normalfont H($\phi$)]: $\phi\colon I\times V\to Y_1$ satisfies
		\begin{description}
			\item[\normalfont (a)] $\phi(t,\cdot)$ is Lipschitz continuous with $L_\phi$, i.e., there exists $L_\phi>0$ satisfying
				$$\|\phi(t,b_1)-\phi(t,b_2)\|_{Y_1}\leq L_\phi\|b_1-b_2\|_V$$
			for all $b_1,\ b_2\in V$ and $t\in I$;
		\end{description}
		\begin{description}
			\item[\normalfont (b)] $\phi(t,0_V)\in L^2(I;Y_1)$.
		\end{description}
	\end{description}

	\begin{description}
		\item[\normalfont H($g$)]: $g\colon Y\times Y\to\mathbb R$ satisfies
		\begin{description}
			\item[\normalfont (a)] $g(\cdot,\cdot)$ is a continuous, bilinear and symmetric form;
		\end{description}
		\begin{description}
			\item[\normalfont (b)] there exist $g_1\in\mathbb R$ and $g_2 >0$ satisfying
				$$g(a,a)+g_1\|a\|^2_{Y_1}\geq g_2\|a\|^2_Y$$
			for all $a\in Y$.
		\end{description}
	\end{description}
	
For the unique solvability of the problem (1.1)-(1.5), first we consider the following sub-problem.

	\noindent\textbf{Problem 3.1.} For the given $w\in C(I;C(\overline\Omega)),\ \xi\in H^1(I;Y_1)\cap L^2(I;Y)$, find $\eta_{w\xi}\colon I\to K_V$ satisfying for each $t\in I$,
	\begin{eqnarray}\label{iq0}
			&&\langle A(\mathcal R(\eta_{w\xi}(t),\xi(t)),\eta_{w\xi}(t)),v-\eta_{w\xi}(t)\rangle_V+j(\mathcal S\eta_{w\xi}(t),\eta_{w\xi}(t),w(t,x),v)\nonumber\\[1mm] 
			&&\qquad\quad-j(\mathcal S\eta_{w\xi}(t),\eta_{w\xi}(t),w(t,x),\eta_{w\xi}(t))\geq\langle f(t),v-\eta_{w\xi}(t)\rangle_V \text{ for all $v\in K_V$}.
	\end{eqnarray}
	
	\begin{theorem}
		\rm Assume H($A$),\ H($j$),\ H($\mathcal R$),\ H($\mathcal S$),\ H($f$) are satisfied and $m>\beta$. Then for any given $w\in C(I;C(\overline\Omega))$ and $\xi\in H^1(I;Y_1)\cap L^2(I;Y)$, \rm{Problem 3.1} has a unique solution $\eta_{w\xi}\in C(I;K_V)$.
	\end{theorem}
	\noindent\textbf{Proof}. For convenience we assume $\mathcal X=C(I;K_V)$ and $\mathcal Y=C(I;K_V)$. For any $u\in C(I;K_V)$, the functions $y_{u\xi}\in C(I;X)$ and $z_u\in C(I;Z)$ are defined as follows: for all $t\in I$,
	\begin{equation*}
		y_{u\xi}(t)=\mathcal R(u(t),\xi(t)),\ z_u(t)=\mathcal Su(t).
	\end{equation*}
	Additionally, the set $\mathcal P$ is defined by
	\begin{eqnarray*}\label{2Phi-u}
	&&\mathcal P= \big\{(\eta_{uw\xi},u)\in\mathcal X\times\mathcal Y \mid \eta_{uw\xi}(t)\in K_V \text{ and} \\
	&&\qquad \quad\langle A(y_{u\xi}(t),\eta_{uw\xi}(t)),v-\eta_{uw\xi}(t)\rangle_V +j(z_u(t),u(t),w(t,x),v)\\
	&&\qquad\qquad-j(z_u(t),u(t),w(t,x),\eta_{uw\xi}(t))\geq\langle f(t),v-\eta_{uw\xi}(t)\rangle_V\ \text{for all }v\in K_V,\ t\in I\big\}.
	\end{eqnarray*}
	Besides, the operator $M:\mathcal X\to\mathcal Y$ is defined as follows:
	\begin{equation*}
		M\eta_{uw\xi}=\eta_{uw\xi}.
	\end{equation*}

Now, we will continue with the following steps.

\smallskip

\textbf{Step 1.} We prove that for any $u\in C(I;K_V)$, there is a unique $\eta_{uw\xi}\in C(I;K_V)$ satisfying
	\begin{eqnarray}\label{iq1}
			&&\langle A(y_{u\xi}(t),\eta_{uw\xi}(t)),v-\eta_{uw\xi}(t)\rangle_V+j(z_u(t),u(t),w(t,x),v)\nonumber\\[1mm]
			&&\qquad -j(z_u(t),u(t),w(t,x),\eta_{uw\xi}(t))\geq\langle f(t),v-\eta_{uw\xi}(t)\rangle_V\text{ for all }v\in K_V,\ t\in I.
	\end{eqnarray}

	In fact, for any fixed $t\in I$, using H($A$) and H($j$)(a), we can conclude that there exists a unique solution $\eta_{uw\xi}(t)$, which can solve (\ref{iq1}) by Lemma 2.2. Now we prove $\eta_{uw\xi}(t)\in C(I;K_V)$. For convenience, considering $t_1,\ t_2\in I$, we write $u(t_i)=u_i$, $\eta_{uw\xi}(t_i)=\eta_i$, $y_{u\xi}(t_i)=y_i$, $z_{u}(t_i)=z_i$, $w(t_i,x)=w_i$, $f(t_i)=f_i$ for $i=1,\ 2$. Using (\ref{iq1}) we get
	\begin{equation}\label{iq2}
		\eta_1\in K_V,\ \langle A(y_1,\eta_1),v-\eta_1\rangle_V+j(z_1,u_1,w_1,v)-j(z_1,u_1,w_1,\eta_1)\geq\langle f_1,v-\eta_1\rangle_V
	\end{equation}
	and
	\begin{equation}\label{iq3}
		\eta_2\in K_V,\ \langle A(y_2,\eta_2),v-\eta_2\rangle_V+j(z_2,u_2,w_2,v)-j(z_2,u_2,w_2,\eta_2)\geq\langle f_2,v-\eta_2\rangle_V
	\end{equation}
	for all $v\in K_V$. Let $v=\eta_2$ and $v=\eta_1$ in (\ref{iq2}) and (\ref{iq3}) respectively. Hence
	\begin{equation}\label{iq4}
		\langle A(y_1,\eta_1),\eta_2-\eta_1\rangle_V+j(z_1,u_1,w_1,\eta_2)-j(z_1,u_1,w_1,\eta_1)\geq\langle f_1,\eta_2-\eta_1\rangle_V
	\end{equation}
	and
	\begin{equation}\label{iq5}
		\langle A(y_2,\eta_2),\eta_1-\eta_2\rangle_V+j(z_2,u_2,w_2,\eta_1)-j(z_2,u_2,w_2,\eta_2)\geq\langle f_2,\eta_1-\eta_2\rangle_V.
	\end{equation}
	In addition, adding (\ref{iq4}) and (\ref{iq5}) we have
	\begin{eqnarray}\label{iq6}
		&&\langle A(y_1,\eta_1)-A(y_2,\eta_2),\eta_2-\eta_1\rangle_V+j(z_1,u_1,w_1,\eta_2) \nonumber\\
		&&\qquad -j(z_1,u_1,w_1,\eta_1)+j(z_2,u_2,w_2,\eta_1)-j(z_2,u_2,w_2,\eta_2)
		\geq\langle f_1-f_2,\eta_2-\eta_1\rangle_V.
	\end{eqnarray}
	Then, it follows from H($A$) that 
	\begin{eqnarray}\label{iq7}
		&&\langle A(y_1,\eta_1)-A(y_2,\eta_2),\eta_2-\eta_1\rangle_V \nonumber\\
		&=&\langle A(y_1,\eta_1)-A(y_1,\eta_2),\eta_2-\eta_1\rangle_V+\langle A(y_1,\eta_2)-A(y_2,\eta_2),\eta_2-\eta_1\rangle_V\nonumber\\
		&\leq&-m\|\eta_1-\eta_2\|^2_V+\|A(y_1,\eta_2)-A(y_2,\eta_2)\|_V \cdot \|\eta_1-\eta_2\|_V \nonumber\\
		&\leq&-m\|\eta_1-\eta_2\|^2_V+L_1\|y_1-y_2\|_X\|\eta_1-\eta_2\|_V.
	\end{eqnarray}
	And using H($j$) we obtain that
	\begin{eqnarray}\label{iq8}
			&&j(z_1,u_1,w_1,\eta_2)-j(z_1,u_1,w_1,\eta_1)+j(z_2,u_2,w_2,\eta_1)-j(z_2,u_2,w_2,\eta_2)\nonumber\\
			&\leq&\alpha\|z_1-z_2\|_Z\|\eta_1-
			\eta_2\|_V+\beta\|u_1-u_2\|_V\|\eta_1-
			\eta_2\|_V+\gamma|w_1-w_2|\|\eta_1-
			\eta_2\|_V.
	\end{eqnarray}
	It follows from (\ref{iq7}), (\ref{iq8}), and (\ref{iq6}) that
	\begin{eqnarray}\label{iq9}
			&&m\|\eta_1-\eta_2\|_V\nonumber\\
			&\leq& L_1\|y_1-y_2\|_X+\alpha\|z_1-z_2\|_Z+\beta\|u_1-u_2\|_V+\gamma|w_1-w_2|-\|f_1-f_2\|_V.
	\end{eqnarray}
	Taking advantage of the definitions of $y_{u\xi}$ and $z_u$, combined with the assumptions of $w$ and H($f$), we show $\eta\in C(I;K_V)$.
	
\textbf{Step 2.} To associate the conclusion of (\ref{iq1}) with Problem 3.1, we need to prove that the operator $\Lambda\colon C(I;K_V)\to C(I;K_V)$ defined by
	\begin{equation*}
		\Lambda u=\eta_{uw\xi}\text{ for all $u\in C(I;V)$}
	\end{equation*}
	has a unique fixed point.
	
	In fact, suppose the set $J\subset I$ is compact. For any $u_1,u_2\in C(I;K_V),\ t\in J$, similar to the derivation of (\ref{iq9}), we get
	\begin{eqnarray}\label{iq10}
			&&m\|\eta_1(t)-\eta_2(t)\|_V\nonumber\\
			&\leq& L_1\|y_1(t)-y_2(t)\|_X+\alpha\|z_1(t)-z_2(t)\|_Z+\beta\|u_1(t)-u_2(t)\|_V.
	\end{eqnarray}
	Now, using H($\mathcal R$) and H($\mathcal S$) we obtain
	\begin{equation}\label{iq11}
		\|y_1(t)-y_2(t)\|_X=\|\mathcal R(u_1,\xi(t))-\mathcal R(u_2,\xi(t))\|_X\leq r_{1J}\int_0^t\|u_1(s)-u_2(s)\|_Vds
	\end{equation}
	and
	\begin{equation}\label{iq12}
		\|z_1(t)-z_2(t)\|_Z=\|\mathcal Su_1(t)-\mathcal Su_2(t)\|_Z\leq s_J\int_0^t\|u_1(s)-u_2(s)\|_Vds.
	\end{equation}
	Then, combined with the inequalities of (\ref{iq11}) and (\ref{iq12}), using the definition of $\Lambda$, from (\ref{iq10}) we have
	\begin{eqnarray*}
			\|\Lambda u_1(t)-\Lambda u_2(t)\|_V&=&\|\eta_1(t)-\eta_2(t)\|_V\nonumber\\
			&\leq& \frac{L_1r_{1J}+\alpha s_J}{m}\int_0^t\|u_1(s)-u_2(s)\|_Vds+\frac{\beta}{m}\|u_1(t)-u_2(t)\|_V.
	\end{eqnarray*}
	
	Combined with the condition that $m>\beta$, $\Lambda$ satisfies Definition 2.2 with $l_K=\frac{\beta}{m}$ and $L_K=\frac{L_1r_{1J}+\alpha s_J}{m}$. So, by Lemma 2.3, $\Lambda$ has a unique fixed point $\eta^*\in C(I;K_V)$ such that $\Lambda\eta^*=\eta^*$.
	
\textbf{Step 3.} Step 1 provides the conditions of Lemma 2.1. Step 2 shows $\Lambda$ has a unique fixed point. So using Lemma 2.1 we know that $\mathcal P\cap Gr(M)$ is a singleton. Therefore, it can be concluded that there exists a unique $\eta^*$ satisfying $(\eta^*,M\eta^*)\in\mathcal P$, so $\eta^*$ solves (\ref{iq0}) and it is unique. That is all the proofs of Theorem 3.1.\hfill$\Box$	
	
	Next, let's consider another sub-problem as follows.
	
	\noindent\textbf{Problem 3.2.} For a given $\xi\in H^1(I;Y_1)\cap L^2(I;Y)$, find $\eta_\xi\colon I\to K_V$, $w_\xi\colon I\times\Omega\to \mathbb R$, such that, the following relationship set up for all $t\in I$:
	
	 \begin{eqnarray}\label{ei0}
	 		\frac{\partial w_\xi(t,x)}{\partial t}-\Delta w_\xi(t,x)=\varphi(\eta_\xi(t),\xi(t)),&&\\
	 		\frac{\partial w_\xi(t,x)}{\partial\nu}\Big|_{\partial\Omega}+bw_\xi(t,x)\Big|_{\partial\Omega}=0,&&\nonumber\\
	 		w_\xi(0,x)=w_0(x),&&\nonumber\\
	 		\langle A(\mathcal R(\eta_\xi(t),\xi(t)),\eta_\xi(t)),v-\eta_\xi(t)\rangle_V+j(\mathcal S\eta_\xi(t),\eta_\xi(t),w_\xi(t,x),v)\quad\qquad&&\nonumber\\
	 		\label{ei1}
	 		-j(\mathcal S\eta_\xi(t),\eta_\xi(t),w_\xi(t,x),\eta_\xi(t))\geq\langle f(t),v-\eta_\xi(t)\rangle_V \text{ for all $v\in K_V$}.&&
	 \end{eqnarray}
	
	\begin{theorem}
		\rm Let $w_0(x)\in C(\overline{\Omega})$ and $m>\beta+T(L_1r_{1J}+\alpha s_J)+\gamma cL_\varphi$. Assume H($\varphi$) and all hypotheses of Theorem 3.1 hold. Then for any given $\xi\in H^1(I;Y_1)\cap L^2(I;Y)$, \rm{Problem 3.2} has a unique solution $(w_\xi,\eta_\xi)\in C(I;C(\overline{\Omega}))\times C(I;K_V)$.
	\end{theorem}
	\noindent\textbf{Proof}. For any given $\xi\in H^1(I;Y_1)\cap L^2(I;Y)$, to solve Problem 3.2, we define $\Theta\colon C(I;K_V)\to C(I;K_V)$ such that $\Theta\eta_\xi(t)=\eta_{w_{\xi\eta}}(t)$, where $\eta_{w_{\xi\eta}}$ solves Problem 3.1 with $w=w_{\xi\eta}$ and $w_{\xi\eta}$ is the solution of equation (\ref{ei0}). Now, we will prove it in accordance with the following steps.
	
\textbf{Step 1.} We prove the unique solvability of equation (\ref{ei0}). 

For any $\eta_\xi\in C(I;K_V)$, let $\bar\varphi(t)=\varphi(\eta_\xi(t),\xi(t))$. Using H($\varphi$)(a) we can get that for all $t\in I$
	\begin{eqnarray}\label{varphi}
			|\bar\varphi(t)|^2&=&|\varphi(\eta_\xi(t),\xi(t))-\varphi(0_V,0_Y)+\varphi(0_V,0_Y)|^2\nonumber\\
			&\leq& 2|\varphi(\eta_\xi(t),\xi(t))-\varphi(0_V,0_Y)|^2+2\varphi(0_V,0_Y)^2\nonumber\\
			&\leq& 2L_\varphi^2\|\eta_\xi(t)\|_V^2+2L_\varphi^2\|\xi(t)\|_{Y_1}^2+2\varphi(0_V,0_Y)^2.
	\end{eqnarray}
	Using H($\varphi$)(b), from (\ref{varphi}) we have
	\begin{eqnarray*}
		&&\|\bar\varphi(t)\|_{L^2(I)}^2\\
		&=&\int_0^T|\varphi(\eta_\xi(t),\xi(t))|^2dt\\
		&\leq&2L_\varphi^2T\|\eta_\xi(t)\|_C^2+2L_\varphi^2\|\xi(t)\|_{L^2(I;Y)}^2+2T\varphi(0_V,0_Y)^2,
	\end{eqnarray*}
	where $\|\eta_\xi(t)\|_C=\sup\limits_{t\in I}\|\eta_\xi(t)\|_V$ is the norm of $C(I;V)$. So $\bar\varphi(t)\in L^2(I)$. In addition, since $\varphi$ is independent of $x$, we can show that $\varphi\in L^a(\Omega)$ and $a$ is arbitrary normal number. Consequently, the conditions of Lemma 2.6 can be satisfied with $a_1=2$ and $b_1=a$. Then, there is a unique $w_{\xi\eta}$ which can solve (\ref{ei0}).

\textbf{Step 2.} We will prove that Problem 3.2 has a unique solution. 

In fact, the proof of Theorem 3.1 shows $\Theta$ is well-defined. To prove the unique solvability of Problem 3.2, we just need to provide the proof for the fixed point of the operator $\Theta$. For convenience, let $w_1=w_{\xi\eta_1}$, $w_2=w_{\xi\eta_2}$ denote the unique solution of (\ref{ei0}) corresponding $\eta_1$ and $\eta_2$ respectively, where we omit $\xi$. Suppose the set $J\subset I$ is compact, so from (\ref{ei1}) we get two inequalities: for any $t\in J$,
	\begin{eqnarray}\label{ei2}
			&&\langle A(\mathcal R(\eta_{w_1}(t),\xi(t)),\eta_{w_1}(t)),v-\eta_{w_1}(t)\rangle_V+j(\mathcal S\eta_{w_1}(t),\eta_{w_1}(t),w_1(t,x),v)\nonumber\\
			&&\quad-j(\mathcal S\eta_{w_1}(t),\eta_{w_1}(t),w_1(t,x),\eta_{w_1}(t))\geq\langle f(t),v-\eta_{w_1}(t)\rangle_V \text{ for all $v\in K_V$}
	\end{eqnarray}
	and
	\begin{eqnarray}\label{ei3}
			&&\langle A(\mathcal R(\eta_{w_2}(t),\xi(t)),\eta_{w_2}(t)),v-\eta_{w_2}(t)\rangle_V+j(\mathcal S\eta_{w_2}(t),\eta_{w_2}(t),w_2(t,x),v)\nonumber\\
			&&\quad-j(\mathcal S\eta_{w_2}(t),\eta_{w_2}(t),w_2(t,x),\eta_{w_2}(t))\geq\langle f(t),v-\eta_{w_2}(t)\rangle_V \text{ for all $v\in K_V$}.
	\end{eqnarray}
	Let $v=\eta_{w_2}$ in (\ref{ei2}) and $v=\eta_{w_1}$ in (\ref{ei3}). Hence
	\begin{eqnarray}\label{ei4}
			&&\langle A(\mathcal R(\eta_{w_1}(t),\xi(t)),\eta_{w_1})-A(\mathcal R(\eta_{w_2}(t),\xi(t)),\eta_{w_2}),\eta_{w_1}(t)-\eta_{w_2}(t)\rangle\nonumber\\
			&&\quad+j(\mathcal S\eta_{w_1}(t),\eta_{w_1}(t),w_1(t,x),\eta_{w_1}(t))-j(\mathcal S\eta_{w_1}(t),\eta_{w_1}(t),w_1(t,x),\eta_{w_2}(t))\nonumber\\
			&&\quad+j(\mathcal S\eta_{w_2}(t),\eta_{w_2}(t),w_2(t,x),\eta_{w_2}(t))-j(\mathcal S\eta_{w_2}(t),\eta_{w_2}(t),w_2(t,x),\eta_{w_1}(t))\leq0.
	\end{eqnarray}
	By the property of inner product, combined with H($A$)(b) and H($\mathcal R$), we have
	\begin{eqnarray}\label{ei5}
			&&\langle A(\mathcal R(\eta_{w_1}(t),\xi(t)),\eta_{w_1})-A(\mathcal R(\eta_{w_2}(t),\xi(t)),\eta_{w_2}),\eta_{w_1}(t)-\eta_{w_2}(t)\rangle\nonumber\\
			&=&\langle A(\mathcal R(\eta_{w_1}(t),\xi(t)),\eta_{w_1})-A(\mathcal R(\eta_{w_1}(t),\xi(t)),\eta_{w_2}),\eta_{w_1}(t)-\eta_{w_2}(t)\rangle\nonumber\\
			&&+\langle A(\mathcal R(\eta_{w_1}(t),\xi(t)),\eta_{w_2})-A(\mathcal R(\eta_{w_2}(t),\xi(t)),\eta_{w_2}),\eta_{w_1}(t)-\eta_{w_2}(t)\rangle\nonumber\\
			&\geq&m\|\eta_{w_1}(t)-\eta_{w_2}(t)\|^2_V-L_1\|\eta_{w_1}(t)-\eta_{w_2}(t)\|_V\|\mathcal R(\eta_{w_1}(t),\xi(t))-\mathcal R(\eta_{w_2}(t),\xi(t))\|_X\nonumber\\
			&\geq&m\|\eta_{w_1}(t)-\eta_{w_2}(t)\|^2_V-L_1r_{1J}\|\eta_{w_1}(t)-\eta_{w_2}(t)\|_V\int_0^t\|\eta_{w_1}(s)-\eta_{w_2}(s)\|_Vds.
	\end{eqnarray}
	Using (\ref{ei5}) and H($j$)(b), it follows from (\ref{ei4}) that
	\begin{eqnarray*}
			&&m\|\eta_{w_1}(t)-\eta_{w_2}(t)\|_V-L_1r_{1J}\int_0^t\|\eta_{w_1}(s)-\eta_{w_2}(s)\|_Vds\nonumber\\
			&\leq&\alpha\|\mathcal S\eta_{w_1}(t)-\mathcal S\eta_{w_2}(t)\|_{L^2}+\beta\|\eta_{w_1}(t)-\eta_{w_2}(t)\|_V+\gamma|w_1(t,x)-w_2(t,x)|.
	\end{eqnarray*}
	Combined with H($\mathcal S$), we conclude that
	\begin{eqnarray}\label{ei7}
			&&\|\eta_{w_1}(t)-\eta_{w_2}(t)\|_V\nonumber\\
			&\leq&\frac{L_1r_{1J}+\alpha s_J}{m-\beta}\int_0^t\|\eta_{w_1}(s)-\eta_{w_2}(s)\|_Vds+\frac{\gamma}{m-\beta}|w_1(t,x)-w_2(t,x)|.
	\end{eqnarray}
	Using Lemma 2.5 we obtain
	\begin{equation}\label{ei8}
		\|w_1-w_2\|_{L^\infty([0,T];L^\infty(\Omega))}\leq c\|\varphi_1-\varphi_2\|_{L^{a_1}([0,T];L^{b_1}(\Omega))}\leq c\|\varphi_1-\varphi_2\|_{L^\infty([0,T];L^\infty(\Omega))},
	\end{equation}
	where $c$ depends on $T,\ N,\ \Omega,\ a$ and the coefficients of the equation (\ref{ei0}).
	
	Since $\varphi$ is continuous, by H($\varphi$) from (\ref{ei8}) we have
	\begin{equation}\label{ei9}
		|w_1(t,x)-w_2(t,x)|\leq cL_{\varphi}\sup\limits_{t\in I}\|\eta_1(t)-\eta_2(t)\|_V=cL_{\varphi}\|\eta_1-\eta_2\|_C.
	\end{equation}

	Next, using (\ref{ei9}) from (\ref{ei7}) we have
	\begin{eqnarray}\label{ei10}
			&&\|\eta_{w_1}(t)-\eta_{w_2}(t)\|_V\nonumber\\
			&\leq&\frac{L_1r_{1J}+\alpha s_J}{m-\beta}\int_0^t\|\eta_{w_1}(s)-\eta_{w_2}(s)\|_Vds+\frac{\gamma cL_{\varphi}}{m-\beta}\|\eta_1-\eta_2\|_C.
	\end{eqnarray}
	Taking the supremum for both sides of the inequality simultaneously, from (\ref{ei10}) and the definition of $\Theta$ we obtain
	\begin{eqnarray}\label{ei11}
			\|\Theta\eta_1-\Theta\eta_2\|_C&=&\|\eta_{w_1}-\eta_{w_2}\|_C\nonumber\\
			&\leq&\frac{L_1r_{1J}+\alpha s_J}{m-\beta}T\|\eta_{w_1}-\eta_{w_2}\|_C+\frac{\gamma cL_{\varphi}}{m-\beta}\|\eta_1-\eta_2\|_C.
	\end{eqnarray}
	Which means 
	\begin{equation*}
		\|\Theta\eta_1-\Theta\eta_2\|_C\leq L\|\eta_1-\eta_2\|_C,
	\end{equation*}
	with $L=\frac{\gamma cL_\varphi}{m-\beta-T(L_1r_{1J}+\alpha s_J)}<1$. By using Banach's fixed point theorem, we can obtain that the operator $\Theta$ has a unique fixed point in $C(I;K_V)$. That is all the proofs of Theorem 3.2.\hfill$\Box$

	Finally, returning to the proof of problem (1.1)-(1.5). For this purpose, we present the relevant results of the parabolic variational inequality (\ref{2}).
	
	\begin{lemma}
		\cite[p.17]{han2016quasistatic} Assume $p\in L^2(I;Y_1)$ and {\normalfont H($g$)} holds. Then there is a unique $\xi\in H^1(I;Y_1)\cap L^2(I;Y)$ satisfying
		\begin{equation}\label{eii3}
			\langle\dot\xi(t),\delta-\xi(t)\rangle_{Y_1}+g(\xi,\delta-\xi)\geq\langle p(t),\delta-\xi(t)\rangle_{Y_1}\text{\rm for all $\delta\in K_Y$},
		\end{equation}
		and $\xi(0)=\xi_0\in K_Y$. Furthermore, for $p=p_i\in L^2(I;Y_1),\ i=1,\ 2$, then we can get the solution $\xi_i$ of {\rm(\ref{eii3})}, and
		\begin{equation*}
			\|\xi_1(t)-\xi_2(t)\|^2_{Y_1}\leq d_1\int_0^t\|p_1(s)-p_2(s)\|^2_{Y_1}ds\ \text{\rm for a.e. }t\in(0,T),
		\end{equation*}
		with a constant $d_1>0$.
	\end{lemma}
	
	Now we will provide the proof of the problem (1.1)-(1.5). For convenience, we denote $\beta+T(L_1r_{1J}+\alpha s_J)+\gamma cL_\varphi+(L_1r_{2J}T+\gamma cL_\varphi)\sqrt{d_1L_\phi T}$ as $m_0$.
	
	\begin{theorem}
		\rm Assume H($A$),\ H($j$),\ H($\mathcal R$),\ H($\mathcal S$),\ H($f$),\ H($\phi$),\ H($g$) and H($\varphi$) are fulfilled. Let $m>m_0$. Then problem (1.1)-(1.5) has a unique solution $(\xi,\eta,w)\in H^1(I;Y_1)\cap L^2(I;Y)\times C(I;K_V)\times C(I;C(\overline{\Omega}))$.
	\end{theorem}
	\noindent\textbf{Proof}. First, let's give the unique solvability of the inequality (\ref{2}). For any given $\kappa\in C(I;K_V)$, by Lemma 3.2, we just need to prove that $\phi\in L^2(I;Y_1)$. Let $\bar\phi(t)=\phi(t,\kappa(t))$. Using H($\phi$)(a) we obtain that for all $t\in I$
	\begin{eqnarray}\label{eii4}
			\|\bar\phi(t)\|^2_{Y_1}&=&\|\phi(t,\kappa(t))-\phi(t,0_V)+\phi(t,0_V)\|^2_{Y_1}\nonumber\\
			&\leq& 2\|\phi(t,\kappa(t))-\phi(t,0_V)\|^2_{Y_1}+2\|\phi(t,0_V)\|^2_{Y_1}\nonumber\\
			&\leq& 2(L_\phi\|\kappa(t)\|_V)^2+2\|\phi(t,0_V)\|^2_{Y_1}\nonumber\\
			&\leq& 2L_\phi^2\|\kappa(t)\|_V^2+2\|\phi(t,0_V)\|^2_{Y_1}.
	\end{eqnarray}
	Using H($\phi$)(b), from (\ref{eii4}) we have
	\begin{equation*}
			\|\bar\phi(t)\|^2_{L^2(I;Y_1)}=\int_0^T\|\phi(t,\kappa(t))\|^2_{Y_1}dt\leq 2L_\phi^2 T\|\kappa(t)\|_C^2+2\|\phi(t,0_V)\|^2_{L^2(I;Y_1)}.
	\end{equation*}
	So $\bar\phi(t)\in L^2(I;Y_1)$.
	
	To solve the Problem 3.1, we define $\Pi\colon C(I;K_V)\to C(I;K_V)$ such that $\Pi\kappa(t)=\eta_{\xi_\kappa,w_{\xi\kappa}}(t)$, where $\xi_\kappa$ solves parabolic variational inequality (\ref{2}) with $\eta=\kappa$, $w_{\xi\kappa}$ is the solution of parabolic equation (\ref{0}) with $\eta=\kappa$ and $\xi=\xi_\kappa$, and $\eta_{\xi_\kappa,w_{\xi\kappa}}$ is the solution of quasivariational inequality (\ref{1}). The proof of Theorems 3.1 and 3.2 shows $\Pi$ is well-defined. To prove the unique solvability of problem (1.1)-(1.5), we just need to provide the proof for the fixed point of the operator $\Pi$. We assume that $\xi_1=\xi_{\kappa_1}, \xi_2=\xi_{\kappa_2}$ and $ w_1=w_{\xi_1\kappa_1}$, $w_2=w_{\xi_2\kappa_2}$ denote respectively the unique solution of (\ref{2}) and (\ref{0}) corresponding $\kappa_1$ and $\kappa_2$ respectively, furthermore we assume that $\eta_1=\eta_{\xi_1,w_1}$ and $\eta_2=\eta_{\xi_2,w_2}$ denote the unique solution of (\ref{1}) corresponding $\kappa_1$ and $\kappa_2$ respectively. Suppose the set $J\subset I$ is compact, in a way analogous to the operator $\Theta$, it can be directly inferred: for all $t\in J$,
	\begin{eqnarray}\label{eii5}
			&&m\|\eta_1(t)-\eta_2(t)\|_V-L_1(r_{1J}\int_0^t\|\eta_1(s)-\eta_2(s)\|_Vds+r_{2J}\int_0^t\|\xi_1(s)-\xi_2(s)\|_{Y_1}ds)\nonumber\\
			&\leq&\alpha\|\mathcal S\eta_1(t)-\mathcal S\eta_2(t)\|_Z+\beta\|\eta_1(t)-\eta_2(t)\|_V+\gamma|w_1(t,x)-w_2(t,x)|.
	\end{eqnarray}
	Considering that H($\mathcal S$) holds, from (\ref{eii5}) it can be derived that
	\begin{eqnarray}\label{eii6}
			&&\|\eta_1(t)-\eta_2(t)\|_V\nonumber\\
			&\leq&\frac{L_1r_{1J}+s_J}{m-\beta}\int_0^t\|\eta_1(s)-\eta_2(s)\|_Vds+\frac{L_1r_{2J}}{m-\beta}\int_0^t\|\xi_1(s)-\xi_2(s)\|_{Y_1}ds\nonumber\\
			&&+\frac{\gamma}{m-\beta}|w_1(t,x)-w_2(t,x)|.
	\end{eqnarray}
	Using Lemma 3.2, considering H($\phi$) we deduce that
	\begin{equation}\label{eii7}
			\|\xi_1(t)-\xi_2(t)\|^2_{Y_1} \leq d_1L_\phi\int_0^t\|\kappa_1(s)-\kappa_2(s)\|^2_Vds.
	\end{equation}
	Similar to (\ref{ei9}), by H($\varphi$) from (\ref{ei8}) we have
	\begin{equation}\label{eii8}
		|w_1(t,x)-w_2(t,x)|\leq cL_{\varphi}\sup\limits_{t\in I}(\|\kappa_1(t)-\kappa_2(t)\|_V+\|\xi_1(t)-\xi_2(t)\|_{Y_1}).
	\end{equation}
	Using (\ref{eii7}) and (\ref{eii8}), from (\ref{eii6}) we have
	\begin{eqnarray}\label{eii9}
			&&\|\eta_1(t)-\eta_2(t)\|_V\nonumber\\
			&\leq&\frac{L_1r_{1J}+s_J}{m-\beta}\int_0^t\|\eta_1(s)-\eta_2(s)\|_Vds+\frac{L_1r_{2J}}{m-\beta}\int_0^t(d_1L_\phi\int_0^s\|\kappa_1(\tau)-\kappa_2(\tau)\|^2_Vd\tau)^\frac{1}{2}ds\nonumber\\
			&&+\frac{\gamma}{m-\beta}cL_{\varphi}\sup\limits_{t\in I}(\|\kappa_1(t)-\kappa_2(t)\|_V+(d_1L_\phi\int_0^t\|\kappa_1(s)-\kappa_2(s)\|^2_Vds)^\frac{1}{2}).
	\end{eqnarray}
	Next, taking the supremum for both sides of the inequality simultaneously, from (\ref{eii9}) we obtain
	\begin{equation*}
			 \|\Pi \kappa_1-\Pi \kappa_2\|_C=\|\eta_1-\eta_2\|_C\leq L\|\kappa_1-\kappa_2\|_C,
	\end{equation*}
	where $L=\frac{\gamma cL_\varphi+(L_1r_{2J}T+\gamma cL_\varphi)\sqrt{d_1L_\phi T}}{m-\beta-T(L_1r_{1J}+s_J)}$. By using Banach's fixed point theorem we can obtain that the operator $\Pi$ has a unique fixed point in $C(I;K_V)$. That is all the proofs of the Theorem 3.3.\hfill$\Box$
\begin{remark}	
	Most PDVIs, the partial differential equations are solved by means of the  semigroup theory. For example, Liu et al. in~\cite{liu2017partial} used the semigroup theory and a fixed point theorem for  a set-valued mapping to establish the existence of solutions for a class of partial differential mixed variational inequalities; Liu et al. in~\cite{Liu2018} applied the semigroup theory and a fixed point theorem of set-valued mapping to show the existence of solutions for a class of partial differential hemivariational inequalities; Li et al. \cite{Li2020} employed operator semigroup theory combined with a fixed point theorem for a set-valued mapping to prove the unique solvability of a class of partial differential set-valued variational inequalities. In contrary, in the paper we prove the unique solvability of the system by using the properties of solution estimates of parabolic partial differential equation and the Banach's fixed point theorem in Sobolev spaces.
\end{remark}

\section{Application to frictional contact problem}\label{S4}

In this section we study the frictional contact problem of viscoelastic materials by applying results of the previous section.

Let $\Omega$ be an open bounded domain in $\mathbb R^d$, $d=2$, $3$ with a smooth boundary $\Gamma$. The boundary $\Gamma$ is partitioned into $\bar\Gamma_1$, $\bar\Gamma_2$ and $\bar\Gamma_3$, where $\Gamma_1$,  $\Gamma_2$ and $\Gamma_3$ are mutually disjoint. Additionally, it is assumed that $\Gamma_1$ satisfies $meas(\Gamma_1)>0$. Vectors and tensors are denoted by boldface letters. For example, the outward unit normal on $\Gamma$ is denoted by $\bm\nu$ and a point in $\mathbb R^d$ is represented as $\bm x=(x_i)$. The time interval is denoted by $I$ which can be $[0,T]$ or 
$\mathbb R_+ = [0,+\infty)$, where $T>0$. Additionally, the derivative of time variable is denoted by the dot above this letter, such as $\dot u$. We define $\mathbb S^d$ as the space of second-order symmetric tensors on $\mathbb R^d$, which is equivalent to the space of $d$-order symmetric matrices. We note that the canonical inner products and the associated norms on $\mathbb R^d$ and $\mathbb S^d$ are defined as
\begin{eqnarray*}
&&\bm u\cdot\bm v=\Sigma u_i v_i,\|\bm u\|=\langle u,u\rangle ^{1/2} \text { for all $\bm u=(u_i),\bm v=(v_i)\in\mathbb R^d$},\\[1mm] &&\bm\sigma\cdot\bm\tau=\Sigma\sigma_i\tau_i,\|\bm\sigma\|=\langle\sigma,\sigma\rangle ^{1/2}\text { for all  $\bm\sigma=(\sigma_{ij}),\bm\tau=(\tau_{ij})\in\mathbb S^d$}.
\end{eqnarray*}
The displacement vector, the stress tensor and the linearized strain tensor are denoted by $\bm u$, $\bm\sigma$ and $\bm\varepsilon(\bm u)$, respectively. The linearized strain tensor is expressed as
\begin{equation*}
\bm\varepsilon_{ij}(\bm u)=\frac{1}{2}(\bm u_{i,j}+\bm u_{j,i}),
\end{equation*}
where $\bm u_{i,j}=\frac{\partial{\bm u_i}}{\partial{\bm x_j}}$. Although, these functions depend on spatial and temporal variables, we do not express this fact  explicitly in the following content that $\bm\sigma$ replaces $\bm\sigma(\bm x)$ or $\bm\sigma(\bm x,t)$.
	
	We study a deformable object, which is on $\Omega$ described above. The entire interior of the object is subject to the body forces of density $\bm f_0(t)$ in $\Omega$ and surface tractions of density $\bm f_2(t)$ in $\Gamma_2$. The object is fixed on $\Gamma_1$, and in contact with a completely rigid body covered with a penetrable soft material of thickness $g$ at the boundary $\Gamma_3$. Over time, this soft material layer will wear out due to friction. More detailed guidelines related to contact friction are provided below.
	
	Now we introduce the mathematical formulation for the following contact problem.

\smallskip
	
	\noindent\textbf{Problem 4.1.} Find the displacement field $\bm u\colon \Omega\times I\to\mathbb R^d$, the stress field $\bm\sigma\colon \Omega\times I\to\mathbb S^d$, the damage field $\xi\colon \Omega\times I\to[0,1]$ and the wear function $w\colon \Gamma_3\times I\to\mathbb R$ 
	such that 
	\begin{eqnarray}
		\bm\sigma(t)=\mathcal{A}\bm\varepsilon(\bm{\dot{u}}(t))+\mathcal{B}\bm\varepsilon(\bm u(t))+\int_{0}^t \mathcal{C}(t-s,\bm\varepsilon(\bm\dot{u}(s)),\xi(s))ds \text{ in } \Omega,&&\\
		\dot{\xi}-\kappa\Delta\xi+\partial I_{[0,1]}(\xi)\ni\phi(\bm\varepsilon(\bm{\dot u}(t))) \text{ in } \Omega,&&\\
		\frac{\partial\xi}{\partial\bm v}=0 \text{ on } \Gamma,&&\\
		Div\bm{\sigma} (t)+\bm{f}_0(t)=\bm{0} \text{ in } \Omega,&&\\
		\bm{u}(t)=\bm{0} \text{ on } \Gamma_1,&&\\
		\bm{\sigma} (t)\bm{\nu}=\bm{f}_2(t) \text{ on } \Gamma_2,&&\\
    	-\bm\sigma_\nu(t)=p(u_\nu(t)-g) \text{ on } \Gamma_3,&&\\
    	\|\bm\sigma_\tau(t)\|\leq\mu(\|\bm{\dot{u}}_\tau(t)\|,w(\bm x,t))|\sigma_\nu(t)|\text{ on } \Gamma_3,&&\\
    	-\bm\sigma_\tau(t)=\mu(\|\bm{\dot{u}}_\tau(t)\|,w(\bm x,t))|\sigma_\nu(t)|\frac{\bm{\dot{u}}_\tau(t)}{\|\bm{\dot{u}}_\tau(t)\|}\text{ on } \Gamma_3,&&\\
    	\frac{\partial w(\bm x,t)}{\partial t}-\Delta w(\bm x,t)=\varphi(\bm\dot u(t),\xi(\bm x,t)) \text{ on } \Gamma_3,&&\\
		\frac{\partial w(\bm x,t)}{\partial\nu}+bw(\bm x,t)=0 \text{ on }\partial\Gamma_3,&&\\
		w(\bm x,0)=w_0(\bm x), {\xi(0)}={\xi_0}\in(0,1) \text{ on } \Gamma_3,&&\\
		\bm u(0)=u_0 \text{ on } \Omega,&&
	\end{eqnarray}
	for each $t\in I$.
	
	(4.1) formulates the viscoelastic constitutive relation, where $\mathcal{A}$ represents the viscosity operator, $\mathcal{B}$ represents the elasticity operator, $\mathcal{C}$ represents the relaxation tensor, and $\bm\varepsilon$ represents the deformation operator. (4.2) formulates the damage function $\xi$ (for details, see \cite{sofonea2005analysis}) , where $\kappa>0$, $\partial I_{[0,1]}$ represents the convex subdifferential of $I_{[0,1]}$, and $\phi$ represents the damage source which depend on $\bm\varepsilon$. Specifically, the damage function can be specified as
	\begin{equation*}
		\phi(\bm\varepsilon(\bm{\dot u}))=\lambda_w-\frac{1}{2}\lambda_E\|\bm\varepsilon(\bm{\dot u})\|^2,
	\end{equation*}
	where $\lambda_E\ and\ \lambda_w$ are positive constant(see \cite{fremond1995damage,fraemond1996damage}). (4.3) is the Neumann condition of $\xi$ on $\Sigma$. (4.4) is the equation of equilibrium without the inertial terms. (4.5)-(4.6) represent the displacement of $0$ on $\Gamma_1$ and traction condition on $\Gamma_2$ respectively. (4.7)-(4.9) formulate the contact with normal compliance and Coulomb's law of dry friction, respectively, and where $\sigma_\nu$ signifies the normal stress, $\bm\sigma_\tau$ denotes the tangential traction, and $\bm{\dot u}_\tau$ represents the tangential component of the velocity field. Additionally, the gap function is denoted by $g$, the normal compliance function is denoted by $p$, and the friction coefficient is denoted by $\mu$ relying on the velocity field and wear. (4.10) describes the wear process, taking into account damage response (for details, see \cite{kalita2019frictional}), $\Delta$ is the Laplace operator $\Gamma_3$, $\varphi$ represents a bilinear functional. This equation is an extended formulation of the Archard law with boundary condition (4.11). Finally, (4.12)-(4.13) represents the initial condition of $w$, $\xi$ and $\bm u$ respectively.
	
	To obtain relevant variational formulations of above problem. The Sobolev space of functions is denoted by $W^{k,p}(\Omega;\mathbb R^d)$ whose components' weak derivative functions up to $k-th$ are $p-th$ integrable on $\Omega$. Furthermore, let $H^k=W^{k,2}(\Omega;\mathbb R^d)$, $H=L^2(\Omega;\mathbb R^d)$ and
	\begin{equation*}
		V=\{\bm v\in H^1(\Omega;\mathbb R^d)|\bm v=0\text{ on }\Gamma_1\},\ Q=\{\bm q=(q_{ij})|\bm q(x)\in\mathbb S^d\text{ and }q_{ij}\in L^2(\Omega;\mathbb R^d)\}\\
	\end{equation*}
	and they respectively have the following inner products
	\begin{equation*}
		\langle \bm v_1,\bm v_2\rangle_V=\int_\Omega\bm\varepsilon(\bm v_1)\cdot\bm\varepsilon(\bm v_2)d\Omega,\ \langle \bm q_1,\bm q_2\rangle_Q=\int_\Omega\bm q_1\cdot\bm q_2 d\Omega.
	\end{equation*}
	In addition, the norms are respectively denoted as $\|\cdot\|_V$ and $\|\cdot\|_Q$, and deformation operator $\bm\varepsilon \colon H^1(\Omega;\mathbb R^d)\to Q$ is denoted by the above definition of the linearized strain tensor. Because of $meas(\Gamma_1)>0$, we can obtain the following relationship by using Korn's inequality,
	\begin{equation}\label{a0}
		\|\bm u\|_V=\|\bm\varepsilon(\bm u)\|_Q \text{ for all }\bm u\in V.
	\end{equation}
	
	For each element $\bm u\in V$, the normal and tangential part are defined as follows
	\begin{equation*}
		u_\nu=\bm u\cdot\bm\nu,\ \bm u_\tau=\bm u-u_\nu\bm\nu.
	\end{equation*}
	Because $u_\nu\in L^2(\Gamma)$, $\bm u_\tau\in L^2(\Gamma;\mathbb R^d)$ and the trace is continuous, it can be derived that
	\begin{equation*}
		\|\bm u\|_{L^2(\Gamma_3;\mathbb R^d)}\leq\|\rho\|\|\bm u\|_V,
	\end{equation*}
	where $\rho\colon V\to L^2(\Gamma_3;\mathbb R^d)$ is the trace operator.
	
	Besides, let $Y=H^1(\Omega;\mathbb R)$ and $Y_1=L^2(\Omega;\mathbb R)$, and the inner products in these two spaces are defined as follows
	\begin{equation*}
		\langle a,b\rangle_{Y_1}=\int_\Omega a\cdot b\ d\Omega,\text{ }\langle c,d \rangle_Y=\langle c,d\rangle_{Y_1}+\int_\Omega \nabla c\cdot\nabla d\ d\Omega.
	\end{equation*}
	In addition, let $K_V=\{u\in V \mid u_\nu\leq g \text{ a.e. on $\Gamma_3$}\}$ and  $K_Y=\{\xi\in Y \mid 0\leq \xi\leq 1 \text{ a.e. on $\Omega$}\}$.
	
	 Next, we start to obtain the variational formulations of the above problem. Suppose $\bm u,\ \bm\sigma,\ \xi$ and $w$ are smooth functions enough and satisfy (4.1)-(4.13). Next, we will use the following Green formula
	\begin{equation*}
		\int_\Omega \bm\sigma\cdot\bm\varepsilon(\bm v)dx+\int_\Omega Div\bm\sigma\cdot\bm vdx=\int_\Gamma\bm\sigma\bm\nu\cdot\bm vd\Gamma \text{  for all $\bm v\in H^1 (\Omega;\mathbb{R}^d)$}.
	\end{equation*}
	Because of the property of $\Gamma$, combined with the Green formula, from the equalities (4.4) and (4.6) we can obtain
	\begin{eqnarray}\label{a1}
			&&\int_\Omega \bm\sigma\cdot(\bm\varepsilon(\bm v)-\bm\varepsilon(\bm{\dot u}(t)))dx=\int_\Omega f_0(t)\cdot(\bm v-\bm{\dot u}(t))dx+\int_{\Gamma_1}\bm\sigma\bm\nu\cdot(\bm v-\bm{\dot u}(t))d\Gamma\nonumber\\
			&&\qquad\quad+\int_{\Gamma_2}\bm f_2(t)\cdot(\bm v-\bm{\dot u}(t))d\Gamma+\int_{\Gamma_3}\bm\sigma\bm\nu\cdot(\bm v-\bm{\dot u}(t))d\Gamma\text{ for all }v\in V.
	\end{eqnarray}
	
	Next, using (4.5) and the conclusion of inner product we can get these identities as follows
	\begin{eqnarray*}
		\bm v-\bm{\dot u}=\bm 0\qquad\qquad\qquad\qquad\qquad\qquad\qquad &&\text{ a.e. on $\Gamma_1$}, \\
		\bm\sigma\bm\nu\cdot(\bm v-\bm{\dot u})=\sigma_\nu(v_\nu-\dot u_\nu)+\bm\sigma_\tau(\bm v_\tau-\bm{\dot u}_\tau) &&\text{ a.e. on $\Gamma_3$}.
	\end{eqnarray*}
	We use these identities to find (\ref{a1}) can be turned into
	\begin{eqnarray}\label{a2}
			&&\int_\Omega \bm\sigma\cdot(\bm\varepsilon(\bm v)-\bm\varepsilon(\bm{\dot u}(t)))dx=\int_\Omega \bm f_0(t)\cdot(\bm v-\bm{\dot u}(t))dx+\int_{\Gamma_2}\bm f_2(t)\cdot(\bm v-\bm{\dot u}(t))d\Gamma\nonumber\\
			&&\qquad\quad+\int_{\Gamma_3}\sigma_\nu\cdot(v_\nu-\dot u_\nu(t))d\Gamma+\int_{\Gamma_3}\bm\sigma_\tau\cdot(\bm v_\tau-\bm{\dot u}_\tau(t))d\Gamma\text{ for all }v\in V.
	\end{eqnarray}
	For $\int_{\Gamma_3}\bm\sigma_\tau\cdot(\bm v_\tau-\bm{\dot u}_\tau(t))d\Gamma$, using (4.8) and (4.9) we obtain
	\begin{eqnarray}\label{a2.1}
		&&\bm\sigma_\tau\cdot(\bm v_\tau-\bm{\dot u}_\tau(t))\nonumber\\
		&=&-\mu(\|\bm{\dot{u}}_\tau(t)\|,w(\bm x,t))|\sigma_\nu|\frac{\bm{\dot{u}}_\tau(t)}{\|\bm{\dot{u}}_\tau(t)\|}\cdot(\bm v_\tau-\bm{\dot u}_\tau(t))\nonumber\\
		&=&-\mu(\|\bm{\dot{u}}_\tau(t)\|,w(\bm x,t))|\sigma_\nu|\frac{\bm{\dot{u}}_\tau(t)}{\|\bm{\dot{u}}_\tau\|}\cdot\bm v_\tau+\mu(\|\bm{\dot{u}}_\tau(t)\|,w(\bm x,t))|\sigma_\nu|\frac{\bm{\dot{u}}_\tau(t)}{\|\bm{\dot{u}}_\tau(t)\|}\cdot\bm{\dot u}_\tau(t)\nonumber\\
		&=&-\mu(\|\bm{\dot{u}}_\tau(t)\|,w(\bm x,t))|\sigma_\nu|\frac{\bm{\dot{u}}_\tau(t)}{\|\bm{\dot{u}}_\tau(t)\|}\cdot\bm v_\tau+\mu(\|\bm{\dot{u}}_\tau(t)\|,w(\bm x,t))|\sigma_\nu|\cdot\|\bm{\dot{u}}_\tau(t)\|\nonumber\\
		&=&\mu(\|\bm{\dot{u}}_\tau(t)\|,w(\bm x,t))|\sigma_\nu|\cdot(\|\bm{\dot{u}}_\tau(t)\|-\frac{\bm{\dot{u}}_\tau(t)\cdot\bm v_\tau}{\|\bm{\dot{u}}_\tau(t)\|})\nonumber\\
		&\geq&\mu(\|\bm{\dot{u}}_\tau(t)\|,w(\bm x,t))|\sigma_\nu|\cdot(\|\bm{\dot{u}}_\tau(t)\|-\|\bm v_\tau\|).
	\end{eqnarray}
	Using (4.1) and (4.7), from (\ref{a2}) and (\ref{a2.1}) we have
	\begin{eqnarray}\label{a3}
			&&\int_\Omega\mathcal{A}\bm\varepsilon(\bm{\dot{u}}(t))\cdot(\bm\varepsilon(\bm v)-\bm\varepsilon(\bm{\dot u}(t)))dx+\int_\Omega\mathcal{B}\bm\varepsilon(\bm u(t))\cdot(\bm\varepsilon(\bm v)-\bm\varepsilon(\bm{\dot u}(t)))dx\nonumber\\
			&&\quad+\int_\Omega\int_{0}^t \mathcal{C}(t-s,\bm\varepsilon(\bm\dot{u}(s)),\xi(s))ds\cdot(\bm\varepsilon(\bm v)-\bm\varepsilon(\bm{\dot u}(t)))dx\nonumber\\
			&&\quad+\int_{\Gamma_3}p(u_\nu(t)-g)(v_\nu-\dot u_\nu(t))d\Gamma\nonumber\\
			&&\quad+\int_{\Gamma_3}\mu(\|\bm{\dot{u}}_\tau(t)\|,w(\bm x,t))p(u_\nu(t)-g)(\|\bm v_\tau\|-\|\bm{\dot{u}}_\tau(t)\|)d\Gamma\nonumber\\
			&\geq&\int_\Omega \bm f_0(t)\cdot(\bm v-\bm{\dot u}(t))dx+\int_{\Gamma_2}\bm f_2(t)\cdot(\bm v-\bm{\dot u}(t))d\Gamma\text{ for all }v\in V.
	\end{eqnarray}
	
	For convenience, the velocity field $\bm{\dot u}$ is denoted by $\bm\eta$. Combined with the initial condition, we obtain
	\begin{equation*}
		\bm u(t)=\mathcal I\bm\eta(t)=\int_0^t \bm\eta(s)ds+\bm u_0 \text{ for all $t\in I$},
	\end{equation*}
	where $\mathcal I\colon C(I;V)\to \mathbb R_+$ is the integral operator.
	With the notation, above variational formulation (\ref{a3}) leads to the following variational formulation of original problem.
	
	\noindent\textbf{Problem 4.2.} Find $\bm\eta\colon I\to V$, $\xi \colon I\to K_Y$ satisfying for each $t\in I$
	\begin{eqnarray}\label{41}
			&&\int_\Omega\mathcal{A}\bm\varepsilon(\bm{\eta}(t))\cdot(\bm\varepsilon(\bm v)-\bm\varepsilon(\bm{\eta}(t)))dx+\int_\Omega\mathcal{B}\bm\varepsilon(\mathcal I\bm\eta(t))\cdot(\bm\varepsilon(\bm v)-\bm\varepsilon(\bm{\bm\eta}(t)))dx\nonumber\\
			&&\quad+\int_\Omega\int_{0}^t \mathcal{C}(t-s,\bm\varepsilon(\bm\eta(s)),\xi(s))ds\cdot(\bm\varepsilon(\bm v)-\bm\varepsilon(\bm{\eta}(t)))dx\nonumber\\
			&&\quad+\int_{\Gamma_3}p((\mathcal I\eta)_\nu(t)-g)(v_\nu-\eta_\nu(t))d\Gamma\nonumber\\
			&&\quad+\int_{\Gamma_3}\mu(\|\bm{\eta}_\tau(t)\|,w(\bm x,t))p((\mathcal I\eta)_\nu(t)-g)(\|\bm v_\tau\|-\|\bm{\eta}_\tau(t)\|)d\Gamma\nonumber\\
			&\geq&\int_\Omega \bm f_0(t)\cdot(\bm v-\bm{\eta}(t))dx+\int_{\Gamma_2}\bm f_2(t)\cdot(\bm v-\bm{\eta}(t))d\Gamma\text{ for all }v\in V.
	\end{eqnarray}
	
	To obtain the variational formulations of Problem 4.1, next we will transform (4.2). Using the properties of convex subdifferential, from (4.2) we have
	\begin{equation}\label{a4}
		\langle \phi(t,\bm\varepsilon(\bm \eta(t)))-\dot\xi(t)+\kappa\Delta\xi,\delta-\xi(t)\rangle_{Y_1} \leq I(\delta)-I(\xi)=0\text{, $\forall\delta\in K_Y$}.
	\end{equation}
	Combining the following formula for integration by parts
	\begin{equation*}
		\int_\Omega(-\Delta u)vd\Omega=\int_\Omega\nabla u\cdot\nabla vd\Omega
	\end{equation*}
	and the definition of bilinear functional $\beta$ denoted by
	\begin{equation}\label{hg}
		g(\delta,\xi)=\kappa\int_\Omega\nabla\delta\cdot\nabla\xi d\Omega,
	\end{equation}
	from (\ref{a4}) we obtain
	\begin{equation}\label{42}
		\langle\dot{\xi(t)},\delta-\xi(t)\rangle_{Y_1}+g(\xi,\delta-\xi)\geq\langle \phi(t,\bm\varepsilon(\bm \eta(t))),\delta-\xi(t)\rangle_{Y_1}\text{, $\forall\delta\in K_Y$}.
	\end{equation}
	
Now we have obtained all variational formulations for the Problem 4.1. For convenience, we write these formulations together and record them as the following problem.
	
	\noindent\textbf{Problem 4.3.} Find $\bm\eta\colon I\to V$, $\xi\colon I\to K_Y$ and $w\colon \Gamma_3\times I\to L^2(\Gamma_3;\mathbb R)$ satisfying for each $t\in I$, the following relationship hold
	\begin{eqnarray}\label{a5}		&&\int_\Omega\mathcal{A}\bm\varepsilon(\bm{\eta}(t))\cdot(\bm\varepsilon(\bm v)-\bm\varepsilon(\bm{\eta}(t)))dx+\int_\Omega\mathcal{B}\bm\varepsilon(\mathcal I\bm\eta(t))\cdot(\bm\varepsilon(\bm v)-\bm\varepsilon(\bm{\bm\eta}(t)))dx\nonumber\\
		&&+\int_\Omega\int_{0}^t \mathcal{C}(t-s,\bm\varepsilon(\bm\eta(s)),\xi(s))ds\cdot(\bm\varepsilon(\bm v)-\bm\varepsilon(\bm{\eta}(t)))dx\nonumber\\
		&&+\int_{\Gamma_3}p((\mathcal I\eta)_\nu(t)-g)(v_\nu-\eta_\nu(t))d\Gamma\nonumber\\
		&&+\int_{\Gamma_3}\mu(\|\bm{\eta}_\tau(t)\|,w(\bm x,t))p((\mathcal I\eta)_\nu(t)-g)(\|\bm v_\tau\|-\|\bm{\eta}_\tau(t)\|)d\Gamma\nonumber\\
		&\geq&\int_\Omega \bm f_0(t)\cdot(\bm v-\bm{\eta}(t))dx+\int_{\Gamma_2}\bm f_2(t)\cdot(\bm v-\bm{\eta}(t))d\Gamma\text{ for all }v\in V.\\		&&\quad\langle\dot{\xi(t)},\delta-\xi(t)\rangle_{Y_1}+g(\xi,\delta-\xi)\geq\langle \phi(t,\bm\varepsilon(\bm \eta(t))),\delta-\xi(t)\rangle_{Y_1}\text{, $\forall\delta\in K_Y$},\nonumber\\
		&&\qquad\qquad\quad\qquad\qquad\qquad\frac{\partial w(\bm x,t)}{\partial t}-\Delta w(\bm x,t)=\varphi(\bm\eta(t),\xi(\bm x,t)) \text{ on } \Gamma_3,\nonumber\\
		&&\qquad\qquad\qquad\qquad\quad\qquad\qquad\qquad w(\bm x,0)=w_0(\bm x), {\xi(0)}={\xi_0}\in(0,1).\nonumber
	\end{eqnarray}
	
To obtain the solvability of Problem 4.3 , we present the following assumptions.

	\textbf{\normalfont H($\mathcal A$)}: $\left\{
	\begin{aligned}
		&\mathcal A\colon \Omega\times\mathbb S^d\to\mathbb S^d\text{ satisfies}\\
		&\text{(a) there exists $L_\mathcal A>0$ satisfying}\\
		&\ \ \ \ \ \ \ \ \ \ \ \ \ \ \ \ \ \ \ \ \|\mathcal A(\bm x,\bm b_1)-\mathcal A(\bm x,\bm b_2)\|\leq L_\mathcal A\|\bm b_1-\bm b_2\|\\
		&\ \ \ \ \ \text{for all }\bm b_1,\bm b_2\in\mathbb S^d,\text{ a.e. }\bm x\in\Omega,\\
		&\text{(b) there exists $m_\mathcal A>0$ satisfying}\\
		&\ \ \ \ \ \ \ \ \ \ \ \ \ \ (\mathcal A(\bm x,\bm b_1)-\mathcal A(\bm x,\bm b_2))\cdot(\bm b_1-\bm b_2)\geq m_\mathcal A\|\bm b_1-\bm b_2\|^2\\
		&\ \ \ \ \ \text{for all } \bm b_1,\bm b_2\in\mathbb S^d,\text{ a.e. }\bm x\in\Omega,\\
		&\text{(c) $\mathcal A(\bm x,\bm b)$ is measurable on $\Omega$, for any $\bm b\in\mathbb S^d$},\\
		&\text{(d) $\mathcal A(\bm x,\bm 0)$ belongs to $Q$}.
	\end{aligned}
	\right.$

	\textbf{\normalfont H($\mathcal B$)}: $\left\{
	\begin{aligned}
		&\mathcal B\colon \Omega\times\mathbb S^d\to\mathbb S^d\text{ satisfies}\\
		&\text{(a) there exists $L_\mathcal B>0$ satisfying}\\
		&\ \ \ \ \ \ \ \ \ \ \ \ \ \ \ \ \ \ \ \ \|\mathcal B(\bm x,\bm b_1)-\mathcal B(\bm x,\bm b_2)\|\leq L_\mathcal B\|\bm b_1-\bm b_2\|\\
		&\ \ \ \ \ \text{for all $\bm b_1,\bm b_2\in\mathbb S^d$, a.e. $\bm x\in\Omega$},\\
		&\text{(b) $\mathcal B(\bm x,\bm b)$ is measurable on $\Omega$, for any $\bm b\in\mathbb S^d$},\\
		&\text{(c) $\mathcal A(\bm x,\bm 0)$ belongs to $Q$}.
	\end{aligned}
	\right.$

	\textbf{\normalfont H($\mathcal C$)}: $\left\{
	\begin{aligned}
		&\mathcal C\colon \Omega\times I\times\mathbb S^d\times\mathbb R\to\mathbb S^d\text{ satisfies}\\
		&\text{(a) there exists $L_\mathcal C>0$ satisfying}\\
		&\ \ \ \ \ \ \ \|\mathcal C(\bm x,t,\bm c_1,d_1)-C(\bm x,t,\bm c_2,d_2)\|\leq L_{\mathcal C}(\|\bm c_1-\bm c_2\|+|d_1-d_2|)\\
		&\ \ \ \ \ \text{for all $\bm c_1,\bm c_2\in\mathbb S^d,\ d_1,d_2\in\mathbb R,\ t\in I$, a.e. $\bm x\in\Omega$},\\
		&\text{(b) $\mathcal C(\bm x,t,\bm c,d)$ is measurable on $\Omega$, for any $t\in I,\bm c\in\mathbb S^d,d\in\mathbb R$},\\
		&\text{(c) $\mathcal C(\bm x,t,\bm 0,d)$ belongs to $Q$}.
	\end{aligned}
	\right.$
	
	\textbf{\normalfont H($p$)}: $\left\{
	\begin{aligned}
		&p\colon \Gamma_3\times\mathbb R\to\mathbb R_+\text{ satisfies}\\
		&\text{(a) there exists $L_p>0$ satisfying}\\
		&\ \ \ \ \ \ \ \ \ \ \ \ \ \ \ \ \ \ \ \ \ \ \ \ |p(\bm x,b_1)-p(\bm x,b_2)|\leq L_p|b_1-b_2|\\
		&\ \ \ \ \ \text{for all $b_1,b_2\in\mathbb R$, a.e. $\bm x\in\Gamma_3$},\\
		&\text{(b) there exists $p^*>0$ satisfying}\\
		&\ \ \ \ p(\bm x,b)\leq p^*\\
		&\ \ \ \ \ \text{for all $b\in\mathbb R$, a.e. $\bm x\in\Gamma_3$},\\
		&\text{(c) $p(\bm x,b)$ is measurable on $\Gamma_3$, for any $b\in\mathbb R$},\\
		&\text{(d) }p(\bm x,b)=0\text{ for all $b\leq 0$, a.e. $\bm x\in\Gamma_3$}.
	\end{aligned}
	\right.$
	
	\textbf{\normalfont H($g'$)}: $
	\begin{aligned}
		&g\in L^2(\Gamma_3)\text{ and }g(\bm x)\geq 0\text{ a.e. }\bm x\in\Gamma_3.
	\end{aligned}
	$
	
	\textbf{\normalfont H($\mu$)}: $\left\{
	\begin{aligned}
		&\mu\colon \Gamma_3\times\mathbb R_+\times\mathbb R\to\mathbb R_+\text{ satisfies}\\
		&\text{(a) there exists $L_\mu>0$ satisfying}\\
		&\ \ \ \ \ \ \ \ \ \ \ \ |\mu(\bm x,b_1,c_1)-\mu(\bm x,b_2,c_2)|\leq L_{\mu}(|b_1-b_2|+|c_1-c_2|)\\
		&\ \ \ \ \ \text{for all $b_1,b_2\in\mathbb R_+,\ c_1,c_2\in\mathbb R$ a.e. $x\in\Gamma_3$},\\
		&\text{(b) there exists $\mu^*>0$ satisfying}\\
		&\ \ \ \ \mu(\bm x,b,c)\leq \mu^*\\
		&\ \ \ \ \ \text{for all $b\in\mathbb R_+,\ c\in\mathbb R$, a.e. $\bm x\in\Gamma_3$},\\
		&\text{(c) $\mu(\bm x,b,c)$ is measurable on $\Gamma_3$, for any $b\in\mathbb R$}.\\
	\end{aligned}
	\right.$
	
	\textbf{\normalfont H($\varphi'$)}: $\left\{
		\begin{aligned}
			&\varphi\colon \Gamma_3\times\mathbb R^d\times\mathbb R\to\mathbb R\text{ satisfies}\\
			&\text{(a) there exists $L\varphi>0$ satisfying}\\
			&\ \ \ \ \ \ \ \ \ |\varphi(\bm x,\bm b_1,c_1)-\varphi(\bm x,\bm b_2,c_2)|\leq \bar L_\varphi(\|\bm b_1-\bm b_2\|+|c_1-c_2|)\\
			&\ \ \ \ \ \text{for all $\bm b_1,\bm b_2\in\mathbb R^d,\ c_1,c_2\in\mathbb R$, a.e. $\bm x\in\Gamma_3$},\\
			&\text{(b) $\varphi(\bm x,\bm b,c)$ is measurable on $\Gamma_3$, for any $\bm b\in\mathbb R^d,\ c\in\mathbb R$}.
		\end{aligned}
	\right.$
	
	\textbf{\normalfont H($\phi'$)}: $\left\{
		\begin{aligned}
			&\phi\colon \Omega\times I\times\mathbb S^d\to\mathbb R\text{ satisfies}\\
			&\text{(a) there exists $\bar L_\phi>0$ satisfying}\\
			&\ \ \ \ \ \ \ \ \ \ \ \ \ \ \ \ \ \ \ |\phi(\bm x,t,\bm c_1)-\phi(\bm x,t,\bm c_2)|\leq \bar L_\phi\|\bm c_1-\bm c_2\|\\
			&\ \ \ \ \ \text{for all $t\in I,\ \bm c_1,\bm c_2\in\mathbb S^d,$ a.e. $\bm x\in\Omega$},\\
			&\text{(b) $\phi(\cdot,\cdot,\bm 0_{\mathbb S^d})\in L^2(I;L^2(\Omega;\mathbb R))$},\\
			&\text{(c) $\phi(\bm x,t,\bm c)$ is measurable on $\Omega$, for any $t\in I,\ \bm c\in\mathbb S^d$}.
		\end{aligned}
	\right.$
	\vskip 0.2cm
	Now, we can provide the proof of the unique solvability of the Problem 4.3.
	\begin{theorem}
		Let the above assumptions {\normalfont H}($\mathcal A$)-{\normalfont H}($\phi'$) hold, $\bm f_2\in C(I;L^2(\Gamma_2;\mathbb R^d))$, and $\bm f_0\in C(I;L^2(\Omega;\mathbb R^d))$. Then the above {\rm Problem 4.3} and {\rm problem (1.1)-(1.5)} have an equivalent relationship.
	\end{theorem}
	
	\noindent\textit{Proof}. For convenience, we define the following operators. Assume $A\colon Q\times V\to V\text{, }\mathcal R\colon C(I;V)\times C(I;Y)\to C(I;Q)\text{, }\mathcal S\colon C(I;V)\to C(I;L^2(\Gamma_3)),\ \phi\colon I\times V\to Y_1$, $j\colon L^2(\Gamma_3)\times V\times\mathbb R\times V\to \mathbb R$ and $f\colon I\to V$ are defined as follows
	\begin{eqnarray}
		\langle A(\bm\theta,\bm u),\bm v\rangle_V&=&\langle\mathcal A\bm\varepsilon(\bm u),\bm\varepsilon(\bm v)\rangle_Q+\langle\bm\theta,\bm\varepsilon(\bm v)\rangle_Q \text{},\label{ha}\\
		\mathcal R(\bm u(t),\xi(t))&=&\mathcal{B}\bm\varepsilon(\mathcal I\bm u(t))+\int_0^t\mathcal{C}(t-s,\bm\varepsilon(\bm u(s)),\xi(s))ds,\label{hr}\\
		\mathcal{S}\bm u(t)&=&(\mathcal I\bm u)_\nu(t)=\int_0^t u_\nu(s)ds+u_{0_\nu},\label{hs}\\
		j(z,\bm u,w,\bm v)&=&\int_{\Gamma_3} p(z-g)v_\nu d\Gamma+\int_{\Gamma_3}\mu(\|\bm u_\tau\|,w)p(z-g)\|\bm v_\tau\|d\Gamma,\label{hj}\\
		\langle f(t),\bm v\rangle_V&=&\int_\Omega\bm f_0(t)\cdot\bm vdx+\int_{\Gamma_2}\bm f_2(t)\cdot\bm vd\Gamma,\label{hf}\\
		\phi(t,\bm\eta(t))&=&\phi(t,\bm\varepsilon(\bm\eta(t))).\label{hphi}
	\end{eqnarray}
	So, (\ref{a5}) can be redescribed finding $\eta\colon I\to V$ satisfying for each $t\in I$
	\begin{eqnarray}\label{a6}
			&&\langle A(\mathcal R(\bm\eta(t),\xi(t)),\bm\eta(t)),\bm v-\bm\eta(t)\rangle_V+j(\mathcal S\bm\eta(t),\bm\eta(t),w(t),\bm v)\nonumber\\
			&&\qquad\quad-j(\mathcal S\bm\eta(t),\bm\eta(t),w(t),\bm\eta(t))\geq\langle f(t),\bm v-\bm\eta(t)\rangle_V \text{ for all $\bm v\in V$}.
	\end{eqnarray}
	
	Note that (\ref{a6}) describe a quasivariational inequality of the form (1.2), where $X=Q$, $V$ and $K_V$ in Section (1)-(3) are equal to $V$ here, and $Z=L^2(\Gamma_3)$. Additionally, $\Omega$ in previous sections is equal to $\Gamma_3$ here. Next, we verify that the operators in Problem 4.1 satisfy the given assumptions in the third section under the above definition.
	
	Using (\ref{ha}), for all $a,\ a_1,\ a_2\in Q$, $b,\ b_1,\ b_2\in V$, we have
	\begin{equation*}
		\langle A(a_1,b)-A(a_2,b),c\rangle_V=\langle a_1-a_2,\bm\varepsilon(c)\rangle_X,
	\end{equation*}
	and
	\begin{equation*}
			\langle A(a,b_1)-A(a,b_2),c\rangle_V=\langle\mathcal A\bm\varepsilon(b_1)-\mathcal A\bm\varepsilon(b_2),\bm\varepsilon(c)\rangle_X.
	\end{equation*}
	Let $c=A(a_1,b)-A(a_2,b)$ and $A(a,b_1)-A(a,b_2)$ respectively, from H($\mathcal A$) we obtain
	\begin{eqnarray*}
		\|A(a_1,b)-A(a_2,b)\|_V^2&=&\langle a_1-a_2,\bm\varepsilon(A(a_1,b)-A(a_2,b))\rangle_X\\
		&\leq&\|a_1-a_2\|_Q\cdot\|\bm\varepsilon(A(a_1,b)-A(a_2,b))\|_X,
	\end{eqnarray*}
	and
	\begin{eqnarray*}
		\|A(a,b_1)-A(a,b_2)\|_V^2&=&\langle\mathcal A\bm\varepsilon(b_1)-\mathcal A\bm\varepsilon(b_2),\bm\varepsilon(A(a,b_1)-A(a,b_2))\rangle_X\\
		&\leq& L_\mathcal A\|\bm\varepsilon(b_1)-\bm\varepsilon(b_2)\|_Q\cdot\|\bm\varepsilon(A(a,b_1)-A(a,b_2))\|_X.
	\end{eqnarray*}
	By the property of $\varepsilon$, so
	\begin{equation*}
		\|A(a_1,b)-A(a_2,b)\|_V\leq\|a_1-a_2\|_X,
	\end{equation*}
	and
	\begin{equation*}
		\|A(a,b_1)-A(a,b_2)\|_V\leq L_\mathcal A\|b_1-b_2\|_V,
	\end{equation*}
	which satisfies the condition H($A$)(a), where $L_1=1,L_2=L_\mathcal A$. Besides, from (\ref{a0}) we obtain
	\begin{eqnarray*}
		\langle A(a,b_1)-A(a,b_2),b_1-b_2\rangle_V&=&\langle\mathcal A\bm\varepsilon(b_1)-\mathcal A\bm\varepsilon(b_2),\bm\varepsilon(b_1-b_2)\rangle_X\\
		&=&\int_\Omega (\mathcal A\bm\varepsilon(b_1)-\mathcal A\bm\varepsilon(b_2))\cdot\bm\varepsilon(b_1-b_2)d\Omega\\
		&\geq& \int_\Omega m_\mathcal A\|\bm\varepsilon(b_1)-\bm\varepsilon(b_2)\|^2d\Omega\\
		&\geq& m_\mathcal A\|b_1-b_2\|^2_V,
	\end{eqnarray*}
	which satisfies the condition H(A)(b), where $m=m_\mathcal A$.
	
	Using (\ref{hj}), from H($p$), H($\mu$) and the property of norm and trace operator, H($j$)(a) is satisfied obviously, in addition, we obtain
	\begin{eqnarray*}
		&&j(a_1,\bm b_1,c_1,\bm d_2)-j(a_1,\bm b_1,c_1,\bm d_1)+j(a_2,\bm b_2,c_2,\bm d_1)-j(a_2,\bm b_2,c_2,\bm d_2)\\
		&=&\int_{\Gamma_3} (p(a_1-g)-p(a_2-g))(\bm d_{2\nu}-\bm d_{1\nu}) d\Gamma\\
		&&+\int_{\Gamma_3}(\mu(\|\bm b_{1\tau}\|,c_1)p(a_1-g)-\mu(\|\bm b_{2\tau}\|,c_2)p(a_2-g))(\|\bm d_{2\tau}\|-\|\bm u_{1\tau}\|)d\Gamma\\
		&\leq&\|p(a_1-g)-p(a_2-g)\|_Z\cdot\|\bm d_{2\nu}-\bm d_{1\nu}\|_Z\\
		&&+\|\mu(\|\bm b_{1\tau}\|,c_1)p(a_1-g)-\mu(\|\bm b_{2\tau}\|,c_2)p(a_2-g)\|_Z\cdot\|(\|\bm d_{2\tau}\|-\|\bm d_{1\tau}\|)\|_Z\\
		&\leq&L_p\|a_1-a_2\|_Z\cdot\|\rho\|\|\bm d_1-\bm d_2\|_V+\mu^*\|p(a_1-g)-p(a_2-g)\|_Z\cdot\|\rho\|\|\bm d_1-\bm d_2\|_V\\
		&&+p^*\|\mu(\|\bm b_{1\tau}\|,c_1)-\mu(\|\bm b_{2\tau}\|,c_2)\cdot\|\rho\|\|\bm d_1-\bm d_2\|_V\\
		&\leq&L_p\|\rho\|\|a_1-a_2\|_Z\|\bm d_1-\bm d_2\|_V+\mu^*L_p\|\rho\|\|a_1-a_2\|_Z\|\bm d_1-\bm d_2\|_V\\
		&&+p^*\|\rho\|(L_{\mu}\|\rho\|\|(\|b_1-b_2\|)\|_V+L_{\mu}|c_1-c_2|)\|\bm d_1-\bm d_2\|_V,
	\end{eqnarray*}
	for all $a_1,\ a_2\in Z,\ \bm{b_1,\ b_2,\ d_1,\ d_2}\in V,\ c_1,\ c_2\in\mathbb R$, which satisfies the condition H($j$)(b), where $\alpha=L_p\|\rho\|(1+\mu^*)$, $\beta=L_{\mu}p^*\|\rho\|^2$ and $\gamma=L_{\mu}p^*\|\rho\|$.
	
	Using (\ref{hr}), let $J\subset I$ and $t\in J$, from H($\mathcal B$) and H($\mathcal C$) we have
	\begin{eqnarray*}
		&&\|\mathcal R(\bm a_1(t),b(t))-\mathcal R(\bm a_2(t),b(t))\|_Q\\
		&\leq&\|\mathcal{B}\bm\varepsilon(\mathcal I\bm a_1(t))-\mathcal B\bm\varepsilon(\mathcal I\bm a_2(t))\|_Q\\
		&&+\int_0^t\|\mathcal{C}(t-s,\bm\varepsilon(\bm a_1(s)),b(s))-\mathcal{C}(t-s,\bm\varepsilon(\bm a_2(s)),b(s))\|_Qds\\
		&\leq&L_\mathcal B\|\mathcal I\bm a_1(t)-\mathcal I\bm a_2(t)\|_V+L_{\mathcal C}\int_0^t\|\bm a_1(s)-\bm a_2(s)\|_Vds\\
		&\leq&(L_{\mathcal B}+L_{\mathcal C})\int_0^t\|\bm a_1(s)-\bm a_2(s)\|_Vds
	\end{eqnarray*}
	and
	\begin{eqnarray*}
		&&\|\mathcal R(\bm a(t),b_1(t))-\mathcal R(\bm a(t),b_2(t))\|_Q\\
		&\leq&\int_0^t\|\mathcal{C}(t-s,\bm\varepsilon(\bm a(s)),b_1(s))-\mathcal{C}(t-s,\bm\varepsilon(\bm a(s)),b_2(s))\|_Qds\\
		&\leq&L_{\mathcal C}\int_0^t\|\bm b_1(s)-\bm b_2(s)\|_Vds,
	\end{eqnarray*}
	for all $\bm{a,\ a_1,\ a_2}\in C(I;V),\ b,\ b_1,\ b_2\in C(I;Y)$, which satisfies the condition H($\mathcal R$), where $r_{1J}=L_{\mathcal B}+L_{\mathcal C}$ and $r_{2J}=L_{\mathcal C}$.
	
	Using (\ref{hs}), let $J\subset I$ and $t\in J$, from the decomposition in normal and tangential part on the boundary we get
	\begin{eqnarray*}
		\|\mathcal S\bm a_1(t)-\mathcal S\bm a_2(t)\|_Z
		&=&\|\int_0^t a_{1\nu}(s)-a_{2\nu}(s)ds\|_Z\\
		&\leq&\int_0^t \|\bm a_1(s)-\bm a_2(s)\|_Vds
	\end{eqnarray*}
	for all $\bm a_1,\ \bm a_2\in C(I;V)$, which satisfies the condition H($\mathcal S$), where $s_J=1$.
	
	Using (\ref{hf}), we can obtain that it satisfies obviously the condition of H($f$).
	
	Using the condition H($\varphi'$) we can obtain directly H($\varphi$), where $L_\varphi=\bar L_\varphi$.
	
	Using (\ref{hphi}), it is obvious that H($\phi$)(b) from H($\phi'$), in addition, combined with inner product in $Y_1$, from (\ref{a0}) we obtain
	\begin{eqnarray*}
		\|\phi(t,\bm b_1)-\phi(t,\bm b_2)\|_{Y_1}^2&=&\int_\Omega(\phi(t,\bm \varepsilon(\bm b_1))-\phi(t,\bm \varepsilon(\bm b_2)))\cdot(\phi(t,\bm b_1)-\phi(t,\bm b_2))d\Omega\\
		&\leq&\bar L_\phi^2\int_\Omega\|\bm \varepsilon(\bm b_1)-\bm \varepsilon(\bm b_2)\|^2d\Omega\\
		&=&\bar L_\phi^2\|\bm b_1-\bm b_2\|_V^2
	\end{eqnarray*}
	for all $\bm b_1,\ \bm b_2\in V,\ t\in I$, which satisfies the condition H($\phi$), where $L_\phi=\bar L_\phi$.

	Finally, by (\ref{hg}) and inner product in $Y$ we obtain
	\begin{equation*}
		\|a\|_Y^2=\|a\|_{Y_1}^2+\int_\Omega\nabla a\cdot\nabla a d\Omega,
	\end{equation*}
	we know that $g(a,a)+\kappa\|a\|_{Y_1}=\kappa\|a\|_Y$ and so H($g$) hold, where $g_1=g_2=\kappa$.

	Now, synthesizing all the above relations and inequalities, we can draw this equivalence relation between Problem (1.1)-(1.5) and Problem 4.3. So by theorem 3.3 we know Problem 4.3 has a unique solution.\hfill$\Box$

\begin{remark}
	Most contact models solved by means of DVIs are formulated for planar contact. For example, Zeng et al. \cite{liu2019differential} investigated a viscoplastic frictionless contact problem; Sofonea et al. \cite{sofonea2016analysis} researched a quasi-static elastic contact problem; Chen et al. \cite{CT2021} studied a elastic frictional contact problem. However, in many practical problems, contact surfaces are curved. Therefore, we investigate a new nonlinear system (1.1)-(1.5) under curved contact and establish the unique solvability of this contact problem by utilizing the theoretical results of this system.
\end{remark}
	
 	\vskip 0.3cm
	\noindent{\bf Data availability}
	
	No data was used for the research described in the article.
	
	\vskip 0.3cm
	\noindent{\bf Declaration of Competing Interest}
	
	We declare that we do not have any commercial or associative interest that represents a conflict of interest in connection
	with the work submitted.
	
	\vskip 0.3cm


	
\end{document}